\documentclass[12pt,reqno]{amsart}

\usepackage{amscd}

\newcommand{\N}{{\mathbb N}}
\newcommand{\Z}{{\mathbb Z}}
\newcommand{\Q}{{\mathbb Q}}
\newcommand{\C}{{\mathbb C}}
\newcommand{\R}{{\mathbb R}}
\renewcommand{\P}{{\mathbb P}}

\newcommand{\EE}{{\mathcal E}}
\newcommand{\FF}{{\mathcal F}}

\newcommand{\Nu}{{\mathcal V}}
\newcommand{\OO}{{\mathcal O}}

\newcommand{\TT}{{\mathcal T}}
\newcommand{\UU}{{\mathcal U}}
\newcommand{\VV}{{\mathcal V}}

\newcommand{\aaa}{{\bf a}}

\newcommand{\rrr}{{\bf r}}
\newcommand{\ddd}{{\rm d}}
\newcommand{\www}{\widetilde}
\newcommand{\oooo}{\overline}

\newcommand{\tm}{{\mathcal T}_M}
\newcommand{\eezz}{\frac{1}{z}}
\newcommand{\nmmm}{{\{0\}\times M}}

\newcommand{\pmmn}{{\P^1\times (M,0)}}
\newcommand{\nmmn}{{\{0\}\times (M,0)}}
\newcommand{\cnnn}{{(\C,0)\times (M,0)}}
\newcommand{\immn}{{\{\infty\}\times (M,0)}}
\newcommand{\cmmn}{{\C\times (M,0)}}
\newcommand{\csmmn}{{\C^*\times (M,0)}}
\newcommand{\csnnn}{{(\C^*,0)\times (M,0)}}

\newcommand{\pmmln}{{\P^1\times (M\times\C^l,0)}}
\newcommand{\nmmln}{{\{0\}\times (M\times\C^l,0)}}
\newcommand{\cnnln}{{(\C,0)\times (M\times\C^l,0)}}
\newcommand{\immln}{{\{\infty\}\times (M\times \C^l,0)}}

\newcommand{\paa}{\partial}

\newcommand{\zdz}{{z\partial_z}}

\newcommand{\nnn}{\nabla}

\DeclareMathOperator{\id}{id}

\DeclareMathOperator{\Lie}{Lie}

\DeclareMathOperator{\rk}{rk}

\DeclareMathOperator{\Res}{Res}

\theoremstyle{plain}
\newtheorem{lemma}{Lemma}[section]
\newtheorem{theorem}[lemma]{Theorem}

\theoremstyle{definition}
\newtheorem{definition}[lemma]{Definition}

\newtheorem{remark}[lemma]{Remark}
\newtheorem{remarks}[lemma]{Remarks}
\newtheorem{example}[lemma]{Example}

\begin{document}

\title[Meromorphic connections and Frobenius manifolds]
{Unfoldings of meromorphic connections and a construction of
Frobenius manifolds} 

\author{
Claus Hertling\and  Yuri Manin}

\address{Claus Hertling,
Max-Planck-Institut f\"ur Mathematik
Vivatsgasse 7, 53111 Bonn, Germany}

\email{hertling\char64 mpim-bonn.mpg.de}

\address{Yuri Manin, Max-Planck-Institut f\"ur Mathematik
Vivatsgasse 7, 53111 Bonn, Germany, 
and Northwestern University, Evanston , USA}

\email{manin@mpim-bonn.mpg.de}



\begin{abstract}
The existence of universal unfoldings of certain germs of meromorphic
connections is established. This is used to prove a 
general construction theorem for Frobenius manifolds. 
A particular case is Dubrovin's theorem on semisimple Frobenius manifolds.
Another special case starts with variations of Hodge structures.
This case is used to compare two constructions of 
Frobenius manifolds, the one in singularity theory and the 
Barannikov--Kontsevich construction. For homogeneous polynomials
which give Calabi--Yau hypersurfaces certain Frobenius submanifolds
in both constructions are isomorphic.
\end{abstract}


\maketitle


\section{Introduction}\label{s1}
\setcounter{equation}{0}

\noindent
Let $\www M$ be a complex manifold. {\it A structure of a Frobenius
manifold on } $\www M$ defined by B. Dubrovin consists of several pieces
of data of which the most important are: 
(a) the choice of a flat structure on $\www M$ represented by a subsheaf
of flat vector fields $\TT^f_{\www M}$ of the sheaf of holomorphic 
vector fields $\TT_{\www M}$; 
(b) a commutative and associative $\OO_{\www M}$-bilinear
multiplication $\circ$ on $\TT_{\www M}$.

Let $\nnn^{(0)}$ be the unique torsionless holomorphic connection
on the holomorphic tangent bundle $T\www M$ with kernel $\TT^f_{\www M}$. 
Then the axiom connecting (a) and (b) is the following requirement:
$\nnn_{z,X} := \nnn_X^{(0)} -\frac{1}{z}X\circ$ is a pencil (with 
parameter $z$) of flat connections on $\www M$, where $X\in \TT_{\www M}$.
It was called the first structure connection in \cite{Man2}.
We use here its more sophisticated version: cf. lemma \ref{t4.4}
and theorem \ref{t4.2}.

Additional pieces of data include a flat holomorphic metric $g$, a flat
identity vector field $e$ for $\circ$, and an Euler field $E$:
see definition \ref{t4.3} below for the complete list of their properties
and interrelations.

Three general constructions of Frobenius manifolds are known: K. Saito's
Frobenius structures on unfolding spaces of singularities;
quantum cohomology (physicists' A-models); and Frobenius structures on the
extended moduli spaces of Calabi--Yau manifolds (Barannikov--Kontsevich
construction, considerably generalized in \cite{Ba1}\cite{Ba2}, 
physicists' B-models).

Some of the important issues of this theory require identification of
Frobenius manifolds obtained by different constructions.
The famous Mirror conjecture of Candelas et al. belongs to this class;
it stimulated much of research.

A general strategy for establishing such an isomorphism consists in 
showing that two Frobenius manifolds to be compared are determined by
a more restricted set of data, and then identifying these data.
For example, if $\www M$ is endowed with a point $0$ at which 
$\circ$ is semisimple and $E\circ $ has a simple spectrum, 
then a finite number of numbers suffices
in order to reconstruct the whole Frobenius germ $(\www M,0)$.
This is a result of Dubrovin \cite[Lecture 3]{Du}.
Generally, however, this is not true; in particular, A- and B-models
fail to satisfy the semisimplicity restriction.

In this case it may happen that for an appropriate submanifold 
$M\subset \www M$ the restriction $H:= T\www M|_M$ endowed with the
restriction of the structure connection and some additional remnants
of Frobenius structure (``initial conditions'') carry sufficient
information in order to reconstruct $\www M$ uniquely.
A statement of this kind will be especially useful if one can ensure that
arbitrary  initial conditions of a given type come from a 
Frobenius structure (``strong reconstruction theorem'').

One of the main results of this paper (theorem \ref{t4.5}) consists
in exhibiting such a set of initial data (see definition \ref{t4.1} (b)).
A special case of this result having direct applications is the theorem
\ref{t5.5} which shows that one can construct a unique Frobenius
manifold from any variation of Hodge filtrations over $M$ and an
opposite filtration assuming that a certain condition which we call
{\it $H^2$-generation} holds (definition \ref{t5.3}).

A prototype of theorem \ref{t5.5} is the first reconstruction theorem
in \cite[Theorem 3.1]{KM}. The uniqueness statement in \ref{t5.5}
is essentially that of \cite{KM}. 
It was already used in \cite[Theorem 6.5]{Ba1}.
The existence is the new part.
The more general theorem \ref{t4.5} is applicable to the germs of 
semisimple Frobenius manifolds $(\www M,0)$ mentioned above as well;
in this case $M$ is the marked point, and one recovers Dubrovin's result.

The initial datum on $M$ can be reformulated as a meromorphic connection
on a bundle on $\P^1\times M$ with a logarithmic pole along 
$\{\infty\}\times M$ and a pole of Poincar\'e rank 1 along
$\{0\}\times M$. So theorem \ref{t4.5} gives a recipe for constructing
Frobenius manifolds: One has to construct such connections.

The proof of theorem \ref{t4.5} relies on the study of unfoldings of
germs of such connections in the chapters \ref{s2} and \ref{s3}.
The main result there is theorem \ref{t2.5}. It shows that a germ of 
a meromorphic connection on a bundle on $\cnnn$ with a pole of
Poincar\'e rank 1 along $\nmmm$ has a universal unfolding if again a certain
generation condition is satisfied. It generalizes the case 
$M=\{pt\}$ which was treated in \cite[(4.1)]{Mal}.

With theorem \ref{t2.5} one can reduce theorem \ref{t4.5} to the case
when $M=\www M$. This case was already known. It was formulated
by Sabbah \cite[Theorem (4.3.6)]{Sab1} \cite[Th\'eor\`eme VII 3.6]{Sab2},
and independently by Barannikov \cite{Ba1}\cite{Ba2}\cite{Ba3}.
A major part of the initial data is equivalent to Barannikov's notion of
{\it a semi-infinite variation of Hodge structures}.
Theorem \ref{t4.5} in the case $M=\www M$ is also implicit in the 
construction in singularity theory \cite{SK}\cite{SM}.

Chapter \ref{s6} is devoted to a class of variations of Hodge structures
to which theorem \ref{t5.5} applies, the variations of Hodge structures
on the primitive parts of
the middle cohomologies of smooth hypersurfaces in $\P^n$ whose degrees
divide $n+1$ (theorem \ref{t6.1}). The corresponding Frobenius manifolds
turn up as submanifolds in two different classes of Frobenius manifolds.

One class arises in singularity theory. The base space of a semiuniversal
unfolding of an isolated hypersurface singularity can be equipped with
the structure of a Frobenius manifold \cite{SK}\cite{SM}\cite{He1}.
Suppose that the singularity is a homogeneous polynomial in 
$\C[x_0,...,x_n]$ of a degree which divides $n+1$.
The submanifold of the base space which parametrizes homogeneous deformations
carries a variation of Hodge structures isomorphic to the one of
the projective hypersurfaces in $\P^n$.
Then a larger submanifold which parametrizes certain semihomogeneous
deformations is a Frobenius submanifold. It is the one determined by
the variation of Hodge filtrations and an opposite filtration.
This is discussed in chapter \ref{s7}.

The other class of Frobenius manifolds arises from the Barannikov--Kontsevich
construction \cite{BK}\cite{Ba1}\cite{Ba2}. There for any Calabi--Yau
manifold a formal germ of a space which extends the moduli space of 
complex structure deformations is constructed and is equipped with a
semi-infinite variation of Hodge structures.
Together with the choice of an opposite filtration this induces the structure
of a Frobenius super manifold on the extended moduli space.
In the case of a Calabi--Yau hypersurface in $\P^n$ it turns out 
\cite[Theorem 6.5]{Ba1} that
the whole semi-infinite variation of Hodge structures is determined by the
variation of Hodge structures on the moduli space of complex
structure deformations.
The whole chapter \ref{s8} is a discussion and reformulation of this result
of Barannikov.

In the case of a homogeneous polynomial of degree $n+1$, the results in 
chapters \ref{s7} and \ref{s8}  give for
suitable opposite filtrations isomorphic Frobenius submanifolds in the
two classes of Frobenius manifolds. 

We thank Dennis Borisov for some discussions about these Frobenius
manifolds and Claude Sabbah for remark \ref{t2.10}.

\bigskip
\noindent
{\bf Index of notations and definitions.}
Both basic structures, that of Frobenius manifolds and initial data,
admit many useful weakenings and variations.
For reader's convenience, we compiled an index of the versions used in 
this paper and related notions.\\
Frobenius manifolds: def. \ref{t4.3}.\\
Frobenius type structure: def. \ref{t4.1} (b).\\
$H^2$-generated germ of a Frobenius manifold of weight $w$: def \ref{t5.4}.\\
$H^2$-generated germ of a variation of filtrations of weight $w$:
def. \ref{t5.3} (a).\\
Higgs field (general): lemma \ref{t2.4}.\\
Higgs field of a Frobenius manifold: def. \ref{t4.3}.\\
$(L)$-structure: def. \ref{t3.1} (c).\\
$(LE)$-structure: remark \ref{t3.3} (i).\\
$(LEP(w))$-structure: def. \ref{t3.1} (b).\\
$(LP(w))$-structure: def. \ref{t3.1} (d).\\
$(T)$-structure: def. \ref{t3.1} (c).\\
$(TE)$-structure: def. \ref{t2.1} (b).\\
$(TEP(w))$-structure: def. \ref{t3.1} (a).\\
$(TP(w))$-structure: def. \ref{t3.1} (d).\\
$(trTLEP(w))$-structure: def. \ref{t4.1} (a).

\section{Unfoldings of meromorphic connections}\label{s2}
\setcounter{equation}{0}

\noindent
Theorem \ref{t2.5} below generalizes a result of Malgrange \cite[(4.1)]{Mal}.
We start with the same setting as in \cite{Mal}.

\begin{definition}\label{t2.1}
(a) Consider a germ $(M,0)$ of a complex manifold, a germ
$H\to \cnnn$ of a holomorphic vector bundle, and a flat connection $\nnn$ 
on the restriction of $H$ to $\csnnn$. Let $z$ be the coordinate on $(\C,0)$. 
The connection $\nnn$ has a pole of Poincar\'e rank $r\in \Z_{\geq 0}$ along
$\nmmn$ if 
\begin{eqnarray}\label{2.1}
\nnn:\OO(H)\to \frac{1}{z^r}\left[ \OO_{\C\times M,0}\cdot \Omega_{M,0}^1 
     + \OO_{\C\times M,0}\cdot \frac{\ddd z}{z}\right] \otimes \OO(H).
\end{eqnarray}
Here $\OO(H)$ is the $\OO_{\C\times M,0}$-module of germs of holomorphic
sections of $H$.
A pole of Poincar\'e rank 0 is called a logarithmic pole.

(b) A $(TE)$-structure is a tuple $((M,0),H,\nnn)$ as in (a) with a pole of
Poincar\'e rank 1.
\end{definition}

\begin{remarks}\label{t2.2}
(i) The results in \cite{Mal} are formulated also for poles of higher 
Poincar\'e rank. We restrict to Poincar\'e rank 1 because that is the case
needed in the later chapters and because a proof of a generalization
of theorem \ref{t2.5}, if possible, would be much more technical.

(ii) The notations $(TE)$, $(TEP(w))$, etc. in the definitions \ref{t2.1}
and \ref{t3.1} are compatible with those in \cite{He2}.
Here '$T$' stands for Twistor, '$E$' for Extension (in $z$-direction),
'$P$' for Pairing, and $w$ is an integer.
\end{remarks}

If $((M,0),H,\nnn)$ is a $(TE)$-structure and $\varphi:(M',0)\to (M,0)$
a holomorphic map of germs of manifolds then one can pull back $H$ and $\nnn$
with $\id\times \varphi:(\C,0)\times (M',0)\to \cnnn$. One easily sees
that $\varphi^*(H,\nnn)$ gives a $(TE)$-structure on $(M',0)$.

\begin{definition}\label{t2.3}
Fix a $(TE)$-structure $((M,0),H,\nnn)$.

(a) An {\it unfolding} of it is a $(TE)$-structure 
$((M\times \C^l,0),\www H,\www
\nnn)$ together with a fixed isomorphism
\begin{eqnarray}\label{2.2}
i:((M,0),H,\nnn)\to ((M\times \C^l,0),\www H, \www\nnn)|_{(M\times\{0\},0)}.
\end{eqnarray}

(b) One unfolding $((M\times \C^l,0),\www H, \www\nnn,i)$ {\it induces} 
a second
unfolding $((M\times \C^{l'},0),\www H', \www\nnn',i')$ if there is an
isomorphism $j$ from the second unfolding to the pullback of the first
unfolding by a map
\begin{eqnarray}\label{2.3}
\varphi:(M\times \C^{l'},0)\to (M\times \C^l,0)
\end{eqnarray}
which is the identity on $(M\times\{0\},0)$, and if
\begin{eqnarray}\label{2.4}
i = j|_{(M\times\{0\},0)}\circ i'.
\end{eqnarray}
(Then $j$ is uniquely determined by $\varphi$ and \eqref{2.4}.)

(c) An unfolding is {\it universal} if it induces any unfolding via a 
unique map $\varphi$.
\end{definition}

By definition, a universal unfolding of a $(TE)$-structure is unique up to
canonical isomorphism if it exists. For existence we need some conditions
which are formulated in terms of the data of the following lemma.

\begin{lemma}\label{t2.4}\cite{Sab2}\cite[ch. 2]{He2}
Let $((M,0),H,\nnn)$ be a $(TE)$-structure. Define the germ of a holomorphic
vector bundle $K:=H|_\nmmn$ on $(M,0)$ with $\OO_{M,0}$-module 
$\OO(K)$ of germs of holomorphic sections. Define maps
\begin{eqnarray}\label{2.5}
C := [z\nnn]&:&\OO(K)\to \Omega^1_{M,0}\otimes\OO(K)\\
&& C_X[a]=[z\nnn_Xa] \mbox{ for } X\in \TT_{M,0},a\in \OO(H)\nonumber
\end{eqnarray}
and
\begin{eqnarray}\label{2.6}
\UU:= [z\nnn_\zdz]:\OO(K)\to \OO(K);
\end{eqnarray}
here $[\ ]$ means restriction to $\nmmn$ and $X\in \TT_{M,0}$ 
is lifted canonically to $\cnnn$.

The maps $C_X$, $X\in \TT_{M,0}$, and $\UU$ are commuting $\OO_{M,0}$-linear 
endomorphisms of $\OO(K)$. Therefore the map $C$ is a Higgs field on $K$.
\end{lemma}

{\it Proof.}
This follows from the flatness of $\nnn$ 
and the fact that $\nnn$ has a pole of Poincar\'e rank 1.

For a general discussion of Higgs fields see e.g. \cite[p. 36]{Sab2}. 
A symmetric Higgs field on the tangent bundle is the same as a
commutative and associative $\OO_M$-bilinear multiplication.
\hfill $\qed$

\bigskip
\begin{theorem}\label{t2.5}
Let $((M,0),H,\nnn)$ be a $(TE)$-structure with data $(K,C,\UU)$ as in lemma
\ref{t2.4}. Suppose that there exists a vector $\zeta\in K_0$ in the fiber 
$K_0$ of $K$ at 0 such that
\begin{list}{}{}
\item[(IC)] 
(injectivity condition) the map
\begin{eqnarray}\label{2.7}
C_\bullet\zeta :T_0M\to K_0,\ \ X\mapsto C_X\zeta
\end{eqnarray}
is injective and 
\item[(GC)] 
(generation condition) $\zeta$ and its images under iteration
of the maps $\UU:K_0\to K_0$ and $C_X:K_0\to K_0$, $X\in T_0M$, generate
$K_0$.
\end{list}
Then a universal unfolding of $((M,0),H,\nnn)$ exists. An unfolding
$((M\times \C^l,0),\www H,\www \nnn,i)$ with data $(\www K,\www C,\www\UU)$
as in lemma \ref{t2.4} is a universal unfolding if and only if the map
\begin{eqnarray}\label{2.8}
\www C_\bullet i(\zeta): T_0(M\times \C^l)\to \www K_0,\ \
X\mapsto \www C_X i(\zeta)
\end{eqnarray}
is an isomorphism.
\end{theorem}

The rest of this chapter is devoted to the proof of this theorem.
In \cite[(4.1)]{Mal} the case $M=\{pt\}$ is considered. In this case
the generation condition $(GC)$ is satisfied if and only if 
$\UU:K_0\to K_0$ has for each eigenvalue only one Jordan block.
For example, if $\UU:K_0\to K_0$ is semisimple, it must have 
simple eigenvalues.
But in chapter \ref{s5} we will use theorem \ref{t2.5} in the case $\UU=0$.

The first part of the proof follows \cite{Mal}, using an extension of $H$
to a bundle on $\pmmn$ and using the rigidity properties of logarithmic
poles, which are formulated in lemma \ref{t2.6}.

\begin{lemma}\label{t2.6}
(a) (e.g. \cite[III.1.20]{Sab2}) 
Let $(L'\to \csnnn,\nnn)$ be the germ of a holomorphic vector bundle with flat
connection $\nnn$, and let $(L^{(0)}\to (\C,0)\times\{0\},\nnn)$ be an
extension of $(L',\nnn)|_{(\C^*,0)\times \{0\}}$ with a logarithmic pole at
0.
Then an extension $(L\to \cnnn,\nnn)$ of $(L',\nnn)$ with a logarithmic
pole along $\{0\}\times(M,0)$ exists with 
$(L,\nnn)|_{(\C,0)\times \{0\}}=(L^{(0)},\nnn)$. It is unique up to 
canonical isomorphism and it is isomorphic to the pullback 
$\varphi^*(L^{(0)},\nnn)$ where $\varphi:(M,0)\to \{0\}$.

(b) (e.g. \cite[5.1]{He2}) Let $(L\to \cnnn,\nnn)$ be (the germ of) a 
holomorphic vector bundle with a flat connection $\nnn$ over 
$\csnnn$ and a logarithmic pole along $\nmmn$.
The bundle $L|_\nmmn$ is equipped with a flat connection $\nnn^{res}$,
the residual connection, and a $\nnn^{res}$-flat endomorphism $\VV^{res}$,
the residue endomorphism.
They are defined by
\begin{eqnarray}\label{2.9}
\nnn^{res}_X[a] &:=& [\nnn_Xa] \mbox{ \ for }X\in \TT_{M,0}, a\in \OO(L),\\
\VV^{res} [a] &:=& [\nnn_\zdz a] \mbox{ \ for }a\in \OO(L).\label{2.10}
\end{eqnarray}
\end{lemma}

\begin{lemma}\label{t2.7}
Let $(H\to \cnnn,\nnn)$ be a germ of a holomorphic vector bundle with a flat
connection $\nnn$ over $\csnnn$. There exists an extension to a trivial 
bundle $H^{(gl)}\to \pmmn$ with flat connection $\nnn$ over
$(\C^*-\{1\})\times (M,0)$, with logarithmic poles along 
$\{1\}\times (M,0)$ and $\immn$, and with trivial monodromy around
$\{1\}\times (M,0)$.
\end{lemma}

{\it Proof.}
We follow \cite[ch. 3]{Mal}. First, one extends $H$ to a bundle
$H^{(\C)}\to \cmmn$ with flat connection $\nnn$ over $\csmmn$. Then one can
choose an extension of $H^{(\C)}$ to a bundle 
$H^{(\P^1)}\to \pmmn$ with a logarithmic pole along $\{\infty\}\times (M,0)$.
This bundle is not necessarily trivial. 

One chooses an $\OO_{M,0}$-basis  $\sigma_1,....,\sigma_{\rk H}$ 
of germs of flat sections of $\OO(H^{(\P^1)})_{(1,0)}$ and defines for
$r=(r_1,...,r_{\rk H})\in \Z^{\rk H}$ a bundle $H^{(r)}\to \pmmn$ as follows:
on $(\P^1-\{1\})\times (M,0)$ it coincides with $H^{(\P^1)}$; the germ
$\OO(H^{(r)})_{(1,0)}$ is generated by the sections $(z-1)^{r_i}\sigma_i$.
The bundle $H^{(r)}$ has a logarithmic pole along $\{1\}\times (M,0)$.
For some $r$ the restriction $H^{(r)}|_{\P^1\times \{0\}}$ is trivial
\cite[(3.2)]{Mal}. Because being a trivial bundle on $\P^1$ is an
open property, the bundle $H^{(gl)}:= H^{(r)}$ for such an $r$ is trivial.
\hfill $\qed$

\begin{remark}\label{t2.8}
If one applies lemma \ref{t2.7} to a $(TE)$-structure $(H\to \cnnn,\nnn)$,
then the extension $(H^{(gl)},\nnn)$ has two distinguished properties:

(i) Because of lemma \ref{t2.6} (a), any unfolding $(\www H\to \cnnln,
\www\nnn)$ of $(H,\www\nnn)$ has a unique extension 
$(\www H^{(gl)}\to \pmmln,\www\nnn)$ with all the properties in lemma 
\ref{2.7} whose restriction to $\P^1\times (M\times \{0\},0)$ coincides
with $(H^{(gl)},\nnn)$.

(ii) Denote by $\nnn^{res}$ the residual connection on $H|_\immn$.
The space 
\begin{eqnarray}\label{2.11}
V(H^{(gl)}) &:=& \{ \mbox{ global hol. sections }v\mbox{ in }H^{(gl)}\\
&&\mbox{ with }\nnn^{res}(v|_\immn)=0\}\nonumber
\end{eqnarray}
is a vector space of dimension $\rk H$. A basis of it is an 
$\OO_{(\P^1,z)\times (M,0)}$-basis of $\OO(H^{(gl)})_{(z,0)}$ for any
$z\in \P^1$. Below we will work with the connection matrix with respect
to such a basis. For an extension $(\www H^{(gl)},\www\nnn)$ as in (i)
the sections in $V(H^{(gl)})$ extend uniquely to the sections in 
$V(\www H^{(gl)})$.
\end{remark}

Using lemma \ref{t2.7} and these observations we can control the unfoldings
of $(TE)$-structures as in theorem \ref{t2.5}. The following lemma is the
main step in its proof.

\begin{lemma}\label{t2.9}
Let $(H^{(gl)}\to \P^1\times (M,0),\nnn)$ be a trivial holomorphic vector
bundle of rank $n\geq 1$ with a flat connection over 
$(\C^*-\{1\})\times (M,0)$, with logarithmic poles along $\{1\}\times (M,0)$
and $\immn$ and with a pole of Poincar\'e rank 1 along 
$\nmmn$. Define $K:=H^{(gl)}|_\nmmn$, $C$ and $\UU$ as in lemma \ref{t2.4}.
Suppose that the generation condition (GC) of theorem \ref{t2.5}
is satisfied for a vector $\zeta\in K_0$. Choose a basis $v_1,...,v_n$
of $V(H^{(gl)})$ with $v_1|_{(0,0)}=\zeta$. Choose $l\in \N$ and 
$n$ functions $f_1,...,f_n\in \OO_{M\times \C^l,0}$ with 
$f_i|_{(M\times\{0\},0)}=0$. Let $(t_1,...,t_m,y_1,...,y_l)=(t,y)$ be 
coordinates on $(M\times \C^l,0)$.

Then there exists a unique unfolding $(\www H^{(gl)}\to \pmmln,\www\nnn)$
of $(H^{(gl)},\nnn)$ with the following properties.
$\www H^{(gl)}$ is a trivial vector bundle with a flat connection over
$(\C^*-\{1\})\times (M\times \C^l,0)$, with logarithmic poles along
$\{1\}\times (M\times \C^l,0)$ and $\immln$ and a pole of Poincar\'e rank
1 along $\nmmln$. Its restriction to $\P^1\times (M\times \{0\},0)$
is $(H^{(gl)},\nnn)$.
Let $\www v_1,...,\www v_n\in V(\www H^{(gl)})$ be the canonical extensions
of $v_1,...,v_n\in V(H^{(gl)})$. Define $\www K := \www H^{(gl)}|_\nmmln$,
$\www C$ and $\www U$ as in lemma \ref{t2.4}. Then
\begin{eqnarray}\label{2.12}
\www C_{\paa/\paa y_\alpha} \www v_1 = \sum_{i=1}^n 
\frac{\paa f_i}{\paa y_\alpha} \www v_i \mbox{ \ for \ }\alpha=1,...,l.
\end{eqnarray}
\end{lemma}

{\it Proof.}
Suppose for a moment that $(\www H^{(gl)},\www\nnn)$ were already constructed.
The connection matrix $\Omega $ with respect to the basis 
$\www v_1,...,\www v_n$,
\begin{eqnarray}\label{2.13}
\www\nnn (\www v_1,...,\www v_n) = (\www v_1,...,\www v_n)\cdot \Omega,
\end{eqnarray}
would take the form
\begin{eqnarray}\label{2.14}
\Omega = \frac{1}{z}\sum_{i=1}^mC_i\ddd t_i + \frac{1}{z}\sum_{\alpha=1}^l
F_\alpha \ddd y_\alpha + (\frac{1}{z^2}U + \frac{1}{z}V + \frac{1}{z-1}W)
\ddd z
\end{eqnarray}
with matrices
\begin{eqnarray}\label{2.15}
C_i,F_\alpha,U,V,W\in M(n\times n,\OO_{M\times \C^l,0}).
\end{eqnarray}
This follows from $\www \nnn^{res}\www v_i|_\immln=0$ and from the pole orders
of $(\www H^{(gl)},\www \nnn)$ along $\{0,1,\infty\}\times (M\times \C^l,0)$.
The flatness condition $\ddd \Omega +\Omega\land \Omega=0$ could be
written as
\begin{eqnarray}\label{2.16}
&& [C_i,C_j]=0\\
&& [C_i,F_\alpha]=0,\label{2.17}\\
&& [F_\alpha,F_\beta]=0,\label{2.18}\\
&& \frac{\paa C_i}{\paa t_j} = \frac{\paa C_j}{\paa t_i},\label{2.19}\\
&& \frac{\paa C_i}{\paa y_\alpha} 
   = \frac{\paa F_\alpha}{\paa t_i},\label{2.20}\\
&& \frac{\paa F_\alpha}{\paa y_\beta} 
   = \frac{\paa F_\beta}{\paa y_\alpha},\label{2.21}
\end{eqnarray}
\begin{eqnarray}
&& [C_i,U]=0,\label{2.22}\\
&& [F_\alpha,U]=0,\label{2.23}\\
&& \frac{\paa U}{\paa t_i} = [V,C_i]-C_i ,\label{2.24}\\
&& \frac{\paa U}{\paa y_\alpha} = [V,F_\alpha]-F_\alpha ,\label{2.25}\\
&& \frac{\paa W}{\paa t_i} = [W,C_i] ,\label{2.26}\\
&& \frac{\paa W}{\paa y_\alpha} = [W,F_\alpha] ,\label{2.27}\\
&& \frac{\paa V}{\paa t_i} = -[W,C_i] ,\label{2.28}\\
&& \frac{\paa V}{\paa y_\alpha} = -[W,F_\alpha] \label{2.29}.
\end{eqnarray}
The condition \eqref{2.12} would mean
\begin{eqnarray}\label{2.30}
(F_\alpha)_{i1} = \frac{\paa f_i}{\paa y_\alpha}.
\end{eqnarray}
The proof consists of three parts. In parts (I) and (II) we restrict 
to the case $l=\alpha=1$. In part (I) we show inductively uniqueness 
and existence of matrices $C_i,F_1,U,V,W$ with 
\eqref{2.16} - \eqref{2.30} and with coefficients
in $\OO_{M,0}[[y_1]]$. In part (II) 
their convergence will be proved with the Cauchy--Kovalevski theorem.
In part (III) the general case will be proved by induction in $l$.

In remark \ref{t2.10} the system of equations 
\eqref{2.16} - \eqref{2.29} will be written in 
a more compact form after some integration.

{\bf Part (I)}. Suppose $l=\alpha=1$ and $y_1=y$.
Define for $w\in \Z_{\geq 0}$
\begin{eqnarray}\label{2.31}
\OO_{M,0}[y]_{\leq w}&:=&\sum_{k=0}^w\OO_{M,0}\cdot y^k,\\
\OO_{M,0}[y]_{>w} &:=& \OO_{M,0}[y]\cdot y^{w+1} ,\label{2.32}\\
\OO_{M,0}[[y]]_{>w} &:=&\OO_{M,0}[[y]]\cdot y^{w+1}\label{2.33}
\end{eqnarray}
and 
\begin{eqnarray}\label{2.34}
M(w)&:=& M(n\times n,\OO_{M,0}\cdot y^w),\\
M(>w)&:=& M(n\times n,\OO_{M,0}[y]_{>w}),\label{2.35}\\ 
M(\leq w) &:=& M(n\times n,\OO_{M,0}[y]_{\leq w}). \label{2.35b}
\end{eqnarray}

{\it Beginning of the induction for $w=0$:}
The connection matrix $\Omega^{(0)}$ of $(H^{(gl)},\nnn)$ with respect to
the basis $v_1,...,v_n$ takes the form
\begin{eqnarray}\label{2.36}
\Omega^{(0)} = \frac{1}{z}\sum_{i=1}^m C_i^{(0)}\ddd t_i + 
\left(\frac{1}{z^2} U^{(0)} + \frac{1}{z} V^{(0)} + \frac{1}{z-1}W^{(0)}\right)
\ddd z
\end{eqnarray}
with matrices $C_i^{(0)}, U^{(0)}, V^{(0)}, W^{(0)}\in M(0)$.
The flatness condition $\ddd \Omega^{(0)} + \Omega^{(0)}\land\Omega^{(0)}=0$
is equivalent to the equations \eqref{2.16}, \eqref{2.19}, \eqref{2.22}, 
\eqref{2.24}, \eqref{2.26}, \eqref{2.28} for the matrices 
$C_i^{(0)},\ U^{(0)},\ V^{(0)},\ W^{(0)}$ instead of 
$C_i,U,V,W$.

{\it Induction hypothesis for $w\in \Z_{\geq 0}$:}
Unique matrices $C_i^{(k)}, U^{(k)}, V^{(k)}, W^{(k)}\in M(k)$ for 
$0\leq k\leq w$ and $F_{1}^{(k)}\in M(k)$ for $0\leq k\leq w-1$
are constructed such that the matrices
\begin{eqnarray}\label{2.37}
C_i^{(\leq w)} := \sum_{k=0}^w C_i^{(k)}\in M(\leq w)
\end{eqnarray}
and the analogously defined matrices 
$U^{(\leq w)}, V^{(\leq w)}, W^{(\leq w)}, F_{1}^{(\leq w-1)}$
satisfy \eqref{2.16}, \eqref{2.19}, \eqref{2.22}, \eqref{2.24},
\eqref{2.26}, \eqref{2.28} modulo $M(>w)$, 
\eqref{2.17}, \eqref{2.20}, \eqref{2.23}, 
\eqref{2.25}, \eqref{2.27}, \eqref{2.29} modulo $M(>w-1)$
and \eqref{2.30} modulo $\OO_{M,0}[[y]]_{>w-1}$.

{\it Induction step from $w$ to $w+1$:}
It consists of three steps:
\begin{list}{}{}
\item[(i)]
Construction of a matrix $F_{1}^{(w)}\in M(w)$ such that the matrix
$F_{1}^{(\leq w)}= F_{1}^{(\leq w-1)}+F_{1}^{(w)}$ together with
the matrices $C_i^{(\leq w)}, U^{(\leq w)}, V^{(\leq w)}, W^{(\leq w)}$
satisfies \eqref{2.17}, \eqref{2.23} modulo $M(>w)$
and \eqref{2.30} modulo $\OO_{M,0}[[y]]_{>w}$.
\item[(ii)]
Construction of matrices 
$C_i^{(w+1)}, U^{(w+1)}, V^{(w+1)}, W^{(w+1)}\in M(w+1)$ such that the 
matrices $C_i^{(\leq w+1)} = C_i^{(\leq w)} + C_i^{(w+1)}$, the analogously
defined matrices $U^{(\leq w+1)}, V^{(\leq w+1)}, W^{(\leq w+1)}$
and the matrix $F_{1}^{(\leq w)}$ satisfy
\eqref{2.20}, \eqref{2.25}, \eqref{2.27}, \eqref{2.29}
modulo $M(>w)$.
\item[(iii)]
Proof of \eqref{2.16}, \eqref{2.19}, \eqref{2.22}, \eqref{2.24}, 
\eqref{2.26}, \eqref{2.28} modulo $M(>w+1)$ for these matrices.
\end{list}

(i) The matrices $C_i^{(\leq w)}$ and $U^{(\leq w)}$ generate an algebra
of commuting matrices in $M(n\times n,\OO_{M,0}[y]/\OO_{M,0}[y]_{>w})$.
Because of the generation condition (GC), the image of the column vector
$(1,0,...,0)^{tr}$ 
under the action of this algebra is the whole 
space $M(n\times 1, \OO_{M,0}[y]/\OO_{M,0}[y]_{>w})$.
This shows two things:
\begin{list}{}{}
\item[$(\alpha)$] 
This algebra contains for any $i=1,...,n$ a matrix $E_i^{(\leq w)}
\in M(n\times n,\OO_{M,0}[y]_{\leq w})$ with first column
\begin{eqnarray}\label{2.38}
(E_i^{(\leq w)})_{j1} =\delta_{ij}.
\end{eqnarray}
\item[$(\beta)$]
Any matrix in $M(n\times n,\OO_{M,0}[y]_{\leq w})$ which commutes with the
matrices $C_i^{(\leq w)}$ and $U^{(\leq w)}$ modulo
$M(>w)$ is modulo $M(>w)$ a linear combination of the matrices
$E_i^{(\leq w)}$ with coefficients in $\OO_{M,0}[y]_{\leq w}$.
\end{list}
Therefore the matrix $F_{1}^{(\leq w)}\in M(\leq w)$ which is defined by
\begin{eqnarray}\label{2.39}
F_{1}^{(\leq w)} =\sum_{i=1}^n \frac{\paa f_i}{\paa y}\cdot
E_i^{(\leq w)} \mbox{ modulo }
M(n\times n,\OO_{M,0}[[y]]_{>w})
\end{eqnarray}
is the unique matrix which satisfies \eqref{2.17}, \eqref{2.23}
modulo $M(>w)$, \eqref{2.30} modulo $\OO_{M,0}[[y]]_{>w}$ and 
$F_{1}^{(\leq w)} = F_{1}^{(\leq w-1)} + F_{1}^{(w)}$ for  some
$F_{1}^{(w)}\in M(w)$.

\medskip
(ii) This step is obvious.

\medskip
(iii) One checks that the derivatives by
$\frac{\paa}{\paa y}$ of the equations
\eqref{2.16}, \eqref{2.19}, \eqref{2.22}, \eqref{2.24}, 
\eqref{2.26}, \eqref{2.28} modulo $M(>w+1)$ hold. For this one uses
\eqref{2.17}, \eqref{2.19}, \eqref{2.20}, \eqref{2.23} - \eqref{2.29} 
modulo $M(>w)$.
For example one calculates modulo $M(>w)$
\begin{eqnarray}\label{2.48}
&& \frac{\paa}{\paa y} [ C_i^{(\leq w+1)},U^{(\leq w+1)}]\\
&\equiv & 
\left[ \frac{\paa F_{1}^{(\leq w)} }{\paa t_i} , U^{(\leq w)}\right]
+ [ C_i^{(\leq w)},[V^{(\leq w)}, F_{1}^{(\leq w)}] - F_{1}^{(\leq w)}]
\nonumber \\
&\equiv &
- \left [ F_{1}^{(\leq w)}, \frac{\paa U^{(\leq w)}}{\paa t_i}\right] 
+ [ F_{1}^{(\leq w)},[V^{(\leq w)}, C_i^{(\leq w)}] - C_i^{(\leq w)}] 
\nonumber \\
&\equiv & 0 \nonumber
\end{eqnarray}
The other calculations are similar or easier. This finishes the proof of 
the induction step from $w$ to $w+1$.
It shows uniqueness and existence of matrices $C_i,F_{1},U,V,W\in
M(n\times n,\OO_{M,0}[[y]])$ with \eqref{2.16} - \eqref{2.30}
and with restrictions
$(C_i,U,V,W)|_{y=0} = (C_i^{(0},U^{(0)}, V^{(0)}, W^{(0)})$.

{\bf Part (II)}. 
We have to show holomorphy of these matrices.
We want to apply the Cauchy--Kovalevski theorem
in the following form (\cite[(1.31), (1.40), (1.41)]{Fo};
there the setting is real analytic, but proofs and statements hold also in the
complex analytic setting):
Given $N\in \N$ and matrices $A_i, B\in M(N\times N, \C\{t_1,...,t_m,y,
x_1,...,x_N\})$ there exists a unique vector
$\Phi\in M(N\times 1,\C\{t_1,...,t_m,y\})$ with
\begin{eqnarray}\label{2.49}
&&\frac{\paa \Phi}{\paa y} = 
\sum_{i=1}^m A_i(t,y,\Phi)\frac{\paa \Phi}{\paa t_i} + B(t,y,\Phi) ,\\
&&\Phi(t,0)= 0.\label{2.50}
\end{eqnarray}
We will construct a system \eqref{2.49} - \eqref{2.50} with $N=(m+3)n^2$
such that it will be satisfied with the entries of the matrices
$C_i-C_i^{(0)}, U-U^{(0)}, V-V^{(0)}, W-W^{(0)}$ as entries of $\Phi$.
The system will be built from the equations \eqref{2.20}, \eqref{2.25},
\eqref{2.27}, \eqref{2.29} and the following equations \eqref{2.51},
\eqref{2.52} with which one can express the entries of $F_1$ as functions
of the entries of $\Phi$.

The commutative subalgebra of $M(n\times n,\OO_{M,0}[[y]])$ which is generated
by the matrices $C_1,...,C_m,U$ is a free $\OO_{M,0}[[y]]$-module of rank 
$n$.
Choose monomials $G^{(j)}$, $j=1,...,n$, in the matrices
$C_1,...,C_m,U$ which form an $\OO_{M,0}[[y]]$-basis of this module.
Then the matrix $(G_{i1}^{(j)})_{ij}$ of the first columns of the 
matrices $G^{(j)}$ is invertible in $M(n\times n,\OO_{M,0}[[y]])$.
Equation \eqref{2.39} gives
\begin{eqnarray}\label{2.51}
F_1 = \sum_{j=1}^n g_j G^{(j)}
\end{eqnarray}
with coefficients $g_j\in \OO_{M,0}[[y]]$ such that
\begin{eqnarray}\label{2.52}
\left( \frac{\paa f_1}{\paa y},...,\frac{\paa f_n}{\paa y}\right)^{tr}
= (G_{i1}^{(j)})\cdot (g_1,...,g_n)^{tr}.
\end{eqnarray}
Replacing the entries of the matrices 
$C_i-C_i^{(0)}, U-U^{(0)}, V-V^{(0)}, W-W^{(0)}$
by indeterminates $x_1,...,x_N$, the coefficients of the matrices $G^{(j)}$ 
become elements of $\C\{t\}[x_1,...,x_N]$, and, 
because of \eqref{2.52}, the coefficients $g_j$ become elements of 
$\C\{t_1,...,t_m,y,x_1,...,x_N\}$. One obtains from \eqref{2.20}, 
\eqref{2.25}, \eqref{2.27}, \eqref{2.29}, \eqref{2.51}, \eqref{2.52}
a system \eqref{2.49} - \eqref{2.50}.
The theorem of Cauchy--Kovalevski shows
$C_i,F_1,U,V,W\in M(n\times n,\OO_{M\times \C,0})$.
This shows lemma \ref{t2.9} in the case $l=1$.

{\bf Part (III)}
By induction in $l$ one obtains  a slightly weaker version of lemma \ref{t2.9},
where \eqref{2.12} is replaced by 
\begin{eqnarray}\label{2.52b}
\left(\www C_{\paa/\paa y_\alpha} 
\www v_1\right)|_{\{y_{\alpha+1}=...=y_l=0\}} 
= \left(\sum_{i=1}^n 
\frac{\paa f_i}{\paa y_\alpha} 
\www v_i \right)|_{\{y_{\alpha+1}=...=y_l=0\}}
\end{eqnarray}
for $\alpha=1,...,l$.
One has a connection matrix $\Omega$ as in \eqref{2.14}. Condition
\eqref{2.52b} is equivalent to \eqref{2.30} with the same restriction to 
$y_{\alpha+1}=...=y_l=0$. But now \eqref{2.21} gives \eqref{2.30}
and \eqref{2.12} for all $y$.
This finishes the proof of lemma \ref{t2.9}.
\hfill $\qed$

\bigskip
{\it Proof of theorem \ref{t2.5}.}
Let $(H\to \cnnn ,\nnn)$ be a $(TE)$-structure with 
$K,C,\UU$ and $\zeta\in K_0$ with all the properties in theorem \ref{t2.5}.
We choose an extension to a trivial bundle with connection
$(H^{(gl)}\to \pmmn,\nnn)$ with the properties in lemma \ref{t2.7},
and we choose a basis $v_1,...,v_n$ of the vector space
$V(H^{(gl)})$ (cf. remark \ref{t2.8}) with $v_1|_{(0,0)}=\zeta$.

By remark \ref{t2.8} (i) an unfolding $(\www H\to \cnnln,\www \nnn,i)$
of $(H,\nnn)$ extends to an unfolding $(\www H^{(gl)}\to
\pmmln,\www\nnn )$ of $(H^{(gl)},\nnn)$ with the properties in lemma
\ref{t2.7}, and the sections $v_1,...,v_n$ extend uniquely to sections
$\www v_1,...,\www v_n$ in the space $V(\www H^{(gl)})$.

Consider the connection matrix $\Omega$ and the matrices 
$C_i,F_\alpha,U,V,W$ as in \eqref{2.13} - \eqref{2.15} for such an unfolding.
Because of \eqref{2.19} - \eqref{2.21} there exists a unique matrix
$A\in M(n\times n,\OO_{M\times \C^l,0})$ with $A(0)=0$ and 
\begin{eqnarray}\label{2.53}
\Omega = \frac{1}{z} \ddd A + (\frac{1}{z^2}U + \frac{1}{z}V + \frac{1}{z-1}W)
\ddd z.
\end{eqnarray}
The matrices $A(t,0),U(t,0),V(t,0),W(t,0)$ are determined by
$(H^{(gl)},\nnn,v_1,...,v_n)$. Lemma \ref{t2.9} says that for an 
{\it arbitrary}
choice of functions $f_i(t,y) = A_{i1}(t,y)-A_{i1}(t,0)$ a {\it unique}
unfolding $(\www H^{(gl)},\www\nnn )$ exists. The first columns of $A$
give a map 
\begin{eqnarray}\label{2.54}
\psi = (A_{11},...,A_{n1}): (M\times \C^l,0)\to (\C^n,0).
\end{eqnarray}
The map in \eqref{2.8} is an isomorphism if and only if $\psi$ is an
isomorphism.

Fix an unfolding $(\www H,\www\nnn)$ such that $\psi$ is an isomorphism.
This is possible thanks to the injectivity condition (IC). Consider
a second unfolding $(\www H'\to (\C,0)\times (M\times \C^{l'},0),\www \nnn')$
of $(H,\nnn)$ with $\www {H'}^{(gl)}, \www v_1',...,\www v_n',\Omega', 
A',\psi'$ defined analogously. If it is induced from
$(\www H,\www\nnn)$ via a map
\begin{eqnarray}\label{2.55}
\varphi : (M\times \C^{l'},0)\to (M\times \C^l,0)
\end{eqnarray}
then $(\id\times \varphi)^* \Omega = \Omega'$ and therefore $A\circ\varphi=A'$
and
\begin{eqnarray}\label{2.56}
\psi\circ\varphi = \psi'.
\end{eqnarray}
This shows that the inducing map $\varphi$ is unique.

If one does not yet know that $(\www H',\www\nnn')$ is induced from 
$(\www H,\www\nnn)$ one can define $\varphi$ by \eqref{2.56} and 
compare the unfoldings $(\www H',\www\nnn')$ and $\varphi^*(\www H,\www\nnn)$.
Since the first columns of the matrix $A'$ and the corresponding matrix
for $\varphi^*(\www H,\www\nnn)$ coincide, lemma \ref{t2.9} says that
the unfoldings are isomorphic. 
This finishes the proof of theorem \ref{t2.5}.

\begin{remark}\label{t2.10}
The system of equations \eqref{2.16} - \eqref{2.29} can be reduced: 
first, by \eqref{2.26} - \eqref{2.29}, $V+W$
is a constant matrix which we denote in this remark by $Res_\infty$. Second,
with $A$ as in \eqref{2.53} one finds for $U$ the formula
\begin{eqnarray}\label{2.58}
U = U(t,0) - (W-W(t,0)) + [Res_\infty, A-A(t,0)]- (A-A(t,0)).
\end{eqnarray}
This formula and the following reduction were shown to us by C. Sabbah.

Now one can transform \eqref{2.16} - \eqref{2.29} to an equivalent system
only in terms of $A,\ W,\ Res_\infty$. We still write it with 
$C_i =\frac{\paa A}{\paa t_i}$, $F_\alpha=\frac{\paa A}{\paa y_\alpha}$, 
$V=Res_\infty-W$ and $U$ given by \eqref{2.58}. 
One finds with some calculations:

{\it The equations \eqref{2.17}, \eqref{2.23}, \eqref{2.27}, and 
the restrictions to $y=0$ of the equations \eqref{2.16}, \eqref{2.22}, 
\eqref{2.24}, \eqref{2.26} are sufficient.} 

The equations \eqref{2.19} - \eqref{2.21} are obvious; 
\eqref{2.16} and \eqref{2.26} (for all $y$)
follow from differentiating them by $\frac{\paa}{\paa y_\alpha}$ and 
some transformations;
the equations \eqref{2.25} and \eqref{2.24} (for all $y$) follow from 
differentiating \eqref{2.58}; 
the equation \eqref{2.22} (for all $y$) follows from 
differentiating it by $\frac{\paa}{\paa y_\alpha}$; 
now \eqref{2.18}
follows with the generation condition (GC) and \eqref{2.28} and 
\eqref{2.29} are obvious.
\end{remark}

\section{Supplements}\label{s3}
\setcounter{equation}{0}

\noindent
For the applications to Frobenius manifolds in chapter \ref{s4} we need
$(TE)$-structures (definition \ref{t2.1}) with an additional ingredient,
a pairing. It is also useful to consider weaker structures, 
$(T)$-structures and $(L)$-structures.
After giving their definitions we will discuss how the concepts and results
of chapter \ref{s2} extend.

\begin{definition}\label{t3.1}\cite[ch. 2]{He2}
(a) Fix $w\in \Z$. A $(TEP(w))$-structure is a $(TE)$-structure 
$(H\to \cnnn,\nnn)$ together with a $\nnn$-flat, $(-1)^w$-symmetric,
nondegenerate pairing
\begin{eqnarray}\label{3.1}
P: H_{(z,t)}\times H_{(-z,t)} \to \C \mbox{ for } 
(z,t)\in (\C^*\times M,0)
\end{eqnarray}
on a representative of $H$ such that the pairing extends to a nondegenerate
$z$-sesquilinear pairing
\begin{eqnarray}\label{3.2}
P:\OO(H)\times \OO(H)\to z^w\OO_{\C\times M,0}.
\end{eqnarray}

(b) An $(LEP(w))$-structure is a germ of a bundle $H\to \cnnn$
with a flat connection $\nnn$ on the restriction to $\csnnn$ 
with a logarithmic pole along $\nmmn$ and with a pairing $P$ with the 
same properties as in (a).

(c) Fix $r\in \Z_{\geq 0}$. Consider a germ of a holomorphic vector bundle
$H\to \cnnn$ with a map
\begin{eqnarray}\label{3.3}
\nnn :\OO(H) \to \frac{1}{z^r}\OO_{\C\times M,0}\cdot \Omega_{M,0}^1
\end{eqnarray}
such that for some representative of $H$ the restrictions 
$H|_{\{z\}\times (M,0)}$, $z\in (\C^*,0)$, are flat connections.
The tuple $((M,0),H,\nnn)$ is called a $(T)$-structure if $r=1$,
and an $(L)$-structure if $r=0$.

(d) A $(T)$-structure with a pairing $P$ with all properties in (a) is 
a $(TP(w))$-structure, an $(L)$-structure with such a pairing is an
$LP(w)$-structure.
\end{definition}

\begin{lemma}\label{t3.2}
Let $(H\to \cnnn,\nnn,P)$ be a $(TEP(w))$-struc\-ture with generation 
condition (GC) (theorem \ref{t2.5}) and let 
$(\www H\to \cnnln,\www\nnn)$ be an unfolding of the underlying
$(TE)$-struc\-ture. Then $P$ extends to $\OO(\www H)$ and 
$(\www H,\www\nnn,P)$ is a $(TEP(w))$-structure.
\end{lemma}

{\it Proof.}
It is sufficient to consider an unfolding in one parameter $y$.
For some representative of $\www H$ the pairing $P$ extends to a
$\nnn$-flat pairing on the restriction to $(\C^*,0)\times (M\times \C,0)$.
We have to show that it takes values on $\OO(\www H)$ in 
$z^w\OO_{\C\times M\times \C,0}$. A priori the values are in 
$\OO_{\C^*\times M\times \C,0}$. Denote $n:=\rk H$ and let $(z,t_1,...,t_m,y)$
be coordinates on $(\C\times M\times \C,0)$. 
Choose an $\OO_{\C\times M\times \C,0}$-basis $\www v_1,...,\www v_n$
of $\OO(\www H)$ with connection matrix
\begin{eqnarray}\label{3.4}
\Omega = \frac{1}{z} \sum_{i=1}^m C_i\ddd t_i + \frac{1}{z}F\ddd y 
+ \frac{1}{z^2}U\ddd z
\end{eqnarray}
with matrices $C_i,F,U\in M(n\times n,\OO_{\C\times M\times \C,0})$,
and the matrix
\begin{eqnarray}\label{3.5}
R := (P(\www v_i,\www v_j)) \in M(n\times n,\OO_{\C^*\times M\times \C,0}).
\end{eqnarray}
Flatness and $z$-sesquilinearity of the pairing give
\begin{eqnarray}\label{3.6}
\ddd R(z,t,y) = \Omega^{tr}(z,t,y)R(z,t,y) + R(z,t,y)\Omega(-z,t,y),
\end{eqnarray}
that means,
\begin{eqnarray}\label{3.7}
z\frac{\paa}{\paa z} R(z,t,y) &=& 
\frac{1}{z} U^{tr}(z,t,y)R(z,t,y) - \frac{1}{z}R(z,t,y)U(-z,t,y),\\
\frac{\paa}{\paa t_i} R(z,t,y) &=& 
\frac{1}{z} C_i^{tr}(z,t,y)R(z,t,y) - \frac{1}{z}R(z,t,y)C_i(-z,t,y),
\label{3.8}\\
\frac{\paa}{\paa y} R(z,t,y) &=& 
\frac{1}{z} F^{tr}(z,t,y)R(z,t,y) - \frac{1}{z}R(z,t,y)F(-z,t,y).\label{3.9}
\end{eqnarray}
Write $R$ as a power series
\begin{eqnarray}\label{3.10}
R(z,t,y) = \sum_{k=0}^\infty R^{(k)}(z,t) \mbox{ with }
R^{(k)}\in M(n\times n,\OO_{\C^*\times M,0}\cdot y^k)
\end{eqnarray}
and define 
\begin{eqnarray}\label{3.11}
R^{(\leq k)} (z,t,y) := \sum_{j=0}^k R^{(j)} (z,t),
\end{eqnarray}
analogously for $C_i,F,U$, with 
$C_i^{(k)},F^{(k)},U^{(k)}\in M(n\times n,\OO_{\C\times M,0}\cdot y^k)$.
Then $R^{(0)}\in M(n\times n,z^w\OO_{\C\times M,0})$ because
$(H,\nnn,P)$ is a $(TEP(w))$-structure.

{\it Induction hypothesis for $k\in \Z_{\geq 0}$:}
\begin{eqnarray}\label{3.12}
R^{(\leq k)} \in M(n\times n,z^w\OO_{\C\times M\times \C,0}).
\end{eqnarray}

{\it Induction step from $k$ to $k+1$:}
Recall the definition of $M(>k)$ in \eqref{2.35}. The equations \eqref{3.7}
and \eqref{3.8} show that one has modulo $M(>k)$
\begin{eqnarray}\label{3.13}
&&C_i^{(\leq k)tr} (0,t,y) [z^{-w}R^{(\leq k)}(z,t,y)]|_{z=0} \nonumber \\
&\equiv & [z^{-w}R^{(\leq k)}(z,t,y)]|_{z=0} C_i^{(\leq k)}(0,t,y) ,\\
&&U^{(\leq k)tr} (0,t,y) [z^{-w}R^{(\leq k)}(z,t,y)]|_{z=0} \nonumber \\
&\equiv & [z^{-w}R^{(\leq k)}(z,t,y)]|_{z=0} U^{(\leq k)}(0,t,y) .\label{3.14}
\end{eqnarray}
Because of the generation condition (GC) the matrix $F^{(\leq k)}(0,t,y)$
is an element of the commutative subalgebra of 
$M(n\times n,\OO_{M,0}[y])/M(>k)$ which is generated by 
$C_1^{(\leq k)},...,C_m^{(\leq k)},U^{(\leq k)}$.
Therefore modulo $M(>k)$
\begin{eqnarray}\label{3.15}
&&F^{(\leq k)tr}(0,t,y) [z^{-w}R^{(\leq k)}(z,t,y)]|_{z=0} \nonumber\\
&\equiv& [z^{-w}R^{(\leq k)}(z,t,y)]|_{z=0} F^{(\leq k)}(0,t,y) .
\end{eqnarray}
This together with \eqref{3.9} completes the induction step.
\hfill $\qed$

\begin{remarks}\label{t3.3}
(i) A bundle $(H\to \cnnn,\nnn)$ with a logarithmic pole along 
$\{0\}\times (M,0)$ is also called an $(LE)$-structure.

(ii) A $(TP(w))$-structure is essentially equivalent to Barannikov's 
notion of a semi-infinite variation of Hodge structures \cite{Ba2}\cite{Ba3}.

(iii) The notion of an unfolding of a $(TE)$-structure (definition \ref{t2.3}
(a)) carries over to structures of type $(TEP(w))$, $(L)$, $(LP(w))$, 
$(LE)$, $(LEP(w))$, $(T)$, $(TP(w))$.

(iv) Lemma \ref{t2.6} (a) says that any unfolding of an $(LE)$-structure
is trivial. The same is true for $(L)$-structures. Therefore the 
analogue of lemma \ref{t3.2} holds for $(LEP(w))$-structures and
$(LP(w))$-structures trivially. An $(L)$-structure comes equipped with a
residual connection $\nnn^{res}$ as in lemma \ref{t2.6}.
In fact, an $(L)$-structure is just a germ of a holomorphic family of flat
connections, parametrized by $(\C,0)$ with the connection 
$\nnn^{res}$ for the parameter $z=0$.

(v) The analogue of lemma \ref{t2.7} for $(T)$-structures is easy:
A $(T)$-structure $(H\to \cnnn,\nnn)$ can be extended to a trivial bundle
$(H^{(gl)}\to \pmmn,\nnn)$ with a holomorphic family of flat connections
on the restrictions 
$H^{(gl)}|_{\{z\}\times (M,0)}$ for $z\in \P^1-\{0\}$. To see this,
one chooses an $\OO_{\C,0}$-basis of sections of 
$\OO(H|_{(\C,0)\times \{0\}})$; one glues $H|_{(\C,0)\times \{0\}}$ to a
trivial bundle on $\P^1-\{0\}$, using this basis; one extends the 
trivial bundle to $(\P^1-\{0\})\times (M,0)$ and glues it to $H$
with $\nnn$.

(vi) A $(T)$-structure $(H\to \cnnn,\nnn)$ comes equipped with a Higgs field
$C$ on $K:= H|_{\{0\}\times (M,0)}$ as in lemma \ref{t2.4}.

(vii) The analogues of theorem \ref{t2.5} and lemma \ref{t2.9} hold for
$(T)$-structures. The proofs are the same, of course without the data
encoding the part of the connection in $z$-direction. The generation condition
reads:\\
(GC)' A vector $\zeta\in K_0$ exists which together with its images
under iterations of the maps $C_X:K_0\to K_0$, $X\in T_0M$, generates $K_0$.

(viii) The analogue of lemma \ref{t3.2} holds for $(T)$-structures with the
generation condition (GC)'. The proof is the same.

(ix) Lemma \ref{t3.4} below holds also for $(TP(w))$-structures and 
$(LP(w))$-structures, of course except \eqref{3.18} and \eqref{3.19}.
\end{remarks}

\begin{lemma}\label{t3.4}
(a) \cite[2.5]{He2} Let $(H\to \cnnn,\nnn,P)$ be a $(TEP(w))$-structure
with $K,C,\UU$ as in lemma \ref{t2.4}. Define a pairing 
$g:\OO(K)\times \OO(K)\to \OO_{M,0}$ by
\begin{eqnarray}\label{3.16}
g([a],[b]) := z^{-w}P(a,b) \mod z\OO_{\C\times M,0} \mbox{ for }a,b\in \OO(H).
\end{eqnarray}
It is $\OO_{M,0}$-bilinear, symmetric, nondegenerate, and it satisfies
\begin{eqnarray}\label{3.17}
g(C_Xa,b)&=&g(a,C_Xb) \mbox{ for }X\in \TT_{M,0},a,b\in \OO(K).\\
g(\UU a,b)&=&g(a,\UU b) \mbox{ for }a,b\in \OO(K).\label{3.18}
\end{eqnarray}

(b) \cite[5.1]{He2} Let $(H\to \cnnn,\nnn,P)$ be an $(LEP(w))$-structure
with $K:= H|_{\{0\}\times (M,0)}$, residual connection $\nnn^{res}$ and 
residue endomorphism $\VV^{res}$ as in lemma \ref{t2.6} (b).  
Define a pairing $g$ as in (a). It is $\OO_{M,0}$-bilinear, symmetric,
nondegenerate, $\nnn^{res}$-flat, and it satisfies
\begin{eqnarray}\label{3.19}
g(\VV^{res}a,b)+ g(a,\VV^{res}b) = w\cdot g(a,b) \mbox{ for }a,b\in \OO(K).
\end{eqnarray}
\end{lemma}

{\it Proof.}
All statements follow easily from the $\nnn$-flatness of $P$ and its other
properties, See \cite[Lemma 2.14 and Lemma 5.3]{He2} for details.
\hfill $\qed$

\section{Construction theorem for Frobenius manifolds}\label{s4}
\setcounter{equation}{0}

\noindent
Associated to a holomorphic Frobenius manifold $\www M$ is a series
of meromorphic connections, parametrized by $w\in\Z$, the {\it (first)
structure connections} (lemma \ref{t4.4}).
Under certain assumptions theorem \ref{t2.5} allows to reconstruct 
anyone of them from its restriction to a submanifold $M\subset\www M$.
This restricted connection is considered as initial datum.
Definition \ref{t4.1} and theorem \ref{t4.2} formalize its properties
in two equivalent ways. Theorem \ref{t4.5} is a construction theorem
for Frobenius manifolds, starting from such an initial datum.
Definition \ref{t4.1} and theorem \ref{t4.2} are also discussed (with
different notations) in \cite[I 1]{Sab1} \cite[VI 2]{Sab2}.

\begin{definition}\label{t4.1}\cite[5.2]{He2}
(a) Fix $w\in \Z$. A $(trTLEP(w))$-structure is a tuple
$((M,0),H,\nnn,P)$ with the following properties.
$(M,0)$ is a germ of a complex manifold; $H\to \pmmn$ is a trivial holomorphic
vector bundle with a flat connection on $H|_\csmmn$; and $P$ is a 
$(-1)^w$-symmetric, nondegenerate, $\nnn$-flat pairing
\begin{eqnarray}\label{4.1}
P: H_{(z,t)}\times H_{(-z,t)}\to \C \mbox{ for }(z,t)\in \csnnn.
\end{eqnarray}
The restriction of $(H,\nnn,P)$ to the germ $\cnnn$ is a 
$(TEP(w))$-structure and the restriction to the germ
$(\P^1,\infty)\times (M,0)$ is an $(LEP(-w))$-structure 
(definition \ref{t3.1}).

(b) A {\it Frobenius type structure} is a tuple $((M,0),K,\nnn^r,C,\UU,\VV,g)$
with the following properties. $(M,0)$ is a germ of a complex manifold;
$K\to (M,0)$ is a germ of a holomorphic vector bundle with flat connection
$\nnn^r$; the map $C:\OO(K)\to \Omega_{M,0}^1\otimes \OO(K)$ is a Higgs 
bundle, i.e., a map such that all the endomorphisms 
$C_X:\OO(K)\to \OO(K)$, $X\in \TT_{M,0}$, commute;
the endomorphism $\UU:\OO(K)\to \OO(K)$ of $K$ satisfies $[C,\UU]=0$;
the endomorphism $\VV:\OO(K)\to \OO(K)$ of $K$ is $\nnn^r$-flat;
and $g:\OO(K)\times \OO(K)\to \OO_{M,0}$ is a symmetric, nondegenerate,
$\nnn^r$-flat pairing. These data satisfy
\begin{eqnarray}\label{4.2}
&& \nnn^r_X(C_Y)-\nnn^r_Y(C_X)-C_{[X,Y]}=0, \\
&& \nnn^r(\UU)-[C,\VV]+C=0, \label{4.3} \\
&& g(C_Xa,b) = g(a,C_Xb),  \label{4.4} \\
&& g(\UU a,b) = g(a,\UU b),  \label{4.5} \\
&& g(\VV a,b) = - g(a,\VV b)  \label{4.6}
\end{eqnarray}
for $X,Y\in \TT_{M,0}, \ a,b\in \OO(K)$.
\end{definition}

\begin{theorem}\label{t4.2}
\cite[I 1]{Sab1}\cite[VI 2]{Sab2}\cite[5.2]{He2}
Fix $w\in \Z$. There is a one-to-one correspondence between 
$(trTLEP(w))$-structures and Frobenius type structures on holomorphic
vector bundles. It is given by the steps in (a) and (b).
They are inverse to one another.

(a) Let $(K \to (M,0),\nnn^r,C,\UU,\VV,g)$ be a Frobenius type structure
on $K$. Let $\pi :\pmmn \to (M,0)$ be the projection. 
Define $H:=\pi^*K$, and let $\psi_z:H_{(z,t)}\to K_t$ for $z\in\P^1$
be the canonical projection.
Extend $\nnn^r,C,\UU,\VV,g$ canonically to $H$. Define
\begin{eqnarray}\label{4.7}
\nnn:= \nnn^r +\eezz C + (\eezz\UU-\Nu+\frac{w}{2}\id )\frac{\ddd z}{z}\ .
\end{eqnarray}
Define a pairing
\begin{eqnarray}\label{4.8}
P:H_{(z,t)}\times H_{(-z,t)}&\to& \C \mbox{ \ \ \ for }(z,t)\in \csnnn\\
(a,b) &\mapsto & z^wg(\psi_z a,\psi_{-z}b)\ .\nonumber
\end{eqnarray}
Then $(H,\nnn,P)$ is a $(trTLEP(w))$-structure.

(b) Let $(H,\nnn,P)$ be a $(trTLEP(w))$-structure. Define $K:=H|_\nmmn$,
$C$ and $\UU$ as in lemma \ref{t2.4} and $g$ as in lemma \ref{3.4} (a).
Let $\nnn^{res}$ and $\VV^{res}$ be the residual connection and the
residue endomorphism on $H|_\immn$ as in lemma \ref{2.6} (b).

Because $H$ is a trivial bundle,
there is a canonical projection $\psi:H\to K$, and the bundles
$K$ and $H|_\immn$ are canonically isomorphic. Structure on 
$H|_\immn$ can be shifted to $K$. Let $\nnn^r$ on $K$ be the shift of
$\nnn^{res}$ and let $\VV$ on $K$ be the shift of 
$\VV^{res}+\frac{w}{2}\id$. 
Then $(K\to M,\nnn^r,C,g,\UU,\VV)$ is a Frobenius type structure and
\eqref{4.7} holds.
\end{theorem}

{\it Proof.}
Part of it follows from the lemmas \ref{t2.4}, \ref{t2.6} (b)
and \ref{t3.4}.
For the rest and for details see \cite[Theorem 5.7]{He2}
or \cite[VI 2]{Sab2}.
\hfill $\qed$

\bigskip
Frobenius type structures and $(trTLEP(w))$-structures can be restricted
to any submanifold of the manifold $(M,0)$ over which they are defined.

\begin{definition}\label{t4.3} (Dubrovin)
A {\it Frobenius manifold} $(M,\circ,e,E,g)$ is a complex manifold $M$
of dimension $\geq 1$ with a commutative and associative 
multiplication $\circ$ on the holomorphic
tangent bundle $TM$ , a unit field $e\in\tm$, 
an {\it Euler field} $E\in \tm$, and
a symmetric nondegenerate $\OO_{M,0}$-bilinear pairing $g$ on $TM$ with the
following properties.
The metric $g$ is multiplication invariant,
\begin{eqnarray}\label{4.9}
g(X\circ Y,Z)=g(Y,X\circ Z)\mbox{ \ for }X,Y,Z\in \tm;
\end{eqnarray}
the Levi--Civita connection $\nnn^g$ of the metric $g$ is flat;
together with the Higgs field 
$C:\tm\to \Omega_{M,0}^1\otimes \tm$ with $C_XY:=-X\circ Y$
it satisfies the potentiality condition
\begin{eqnarray}\label{4.10}
\nnn^g_X (C_Y) - \nnn^g_Y(C_X) - C_{[X,Y]}=0 \mbox{ \ for }X,Y\in \tm;
\end{eqnarray}
the unit field $e$ is $\nnn^g$-flat;
the Euler field satisfies $\Lie_E(\circ)=\circ$ and 
$\Lie_E(g)=(2-d)\cdot g$ for some $d\in\C$.
\end{definition}

\begin{lemma}\label{t4.4}
(Structure connections of a Frobenius manifold, e.g. 
\cite[Lecture 3]{Du}, \cite[I 2.5.2]{Man2}, 
\cite[VII 1]{Sab2}, \cite[Lemma 5.11]{He2})
Let $((M,0),\circ,e,E,g)$ be the germ of a Frobenius manifold with Higgs field
$C$, Levi-Civita connection $\nnn^g$ and $d\in\C$ 
as in definition \ref{t4.3}.
Define the endomorphisms $\UU:=E\circ:\TT_{M,0}\to \TT_{M,0}$ and 
\begin{eqnarray}\label{4.11}
\VV:= \TT_{M,0}\to \TT_{M,0}, \ X\mapsto \nnn^g_X E - \frac{2-d}{2}X.
\end{eqnarray}
Then $(TM,\nnn^g,C,\UU,\VV,g)$ is a Frobenius type structure on $TM$.
The unit field $e$ satisfies $\nnn^g e=0$ and $\VV e=\frac{d}{2}e$.

For any $w\in \Z$ theorem \ref{t4.2} gives a $(trTLEP(w))$-structure
$((M,0),H=\pi^*TM,\nnn,P)$, where $\pi :\pmmn\to (M,0)$ is the projection.
These structures are called (first) structure connections of the
Frobenius manifold.
\end{lemma}

One can recover a Frobenius manifold from a structure connection and
the unit field. More abstractly, one can construct from a
$(trTLEP(w))$-structure and a global section with sufficiently nice
properties a unique Frobenius manifold such that the $(trTLEP(w))$-structure
and the global section are isomorphic to a structure connection and the
unit field \cite{Sab1}\cite{Sab2}\cite{Ba1}\cite{Ba2}\cite{Ba3}
(cf. remark \ref{t4.6} (i)). Theorem \ref{t2.5}
allows to start under certain assumptions with a $(trTLEP(w))$-structure
and a global section over a {\it smaller} base space, 
to unfold them universally and then get a Frobenius manifold.
This is formulated in theorem \ref{t4.5} in terms of Frobenius type structures.

\begin{theorem}\label{t4.5}
(Construction theorem for Frobenius manifolds)
Let $((M,0),K,\nnn^r,C,\UU,\VV,g)$ be a Frobenius type structure and $\zeta
\in K_0$ a vector with the following properties:
\begin{list}{}{}
\item[(IC)] (injectivity condition) the map 
$C_\bullet\zeta:T_0M\to K_0$, $X\mapsto C_X\zeta$ is injective.
\item[(GC)] (generation condition) $\zeta$ and its images under iteration
of the maps $C_X:K_0\to K_0$, $X\in T_0M$, and $\UU:K_0\to K_0$ generate
$K_0$.
\item[(EC)] (eigenvector condition) $\VV\zeta=\frac{d}{2}\zeta$ for some
$d\in \C$.
\end{list}

Then there exist up to canonical isomorphism unique data
$((\www M,0),\circ,e,E,\www g,i,j)$ with the following properties.
$((\www M,0),\circ,e,E,\www g)$ is a germ of a Frobenius manifold,
$i:(M,0)\to (\www M,0)$ is an embedding, $j:K\to T\www M|_{i(M)}$ is an
isomorphism above $i$ of germs of vector bundles which maps
$\zeta$ to $e|_0$ and which identifies the Frobenius type structure on $K$
with the natural Frobenius type structure on $T\www M|_{i(M)}$ which
is induced by that on $T\www M$.
\end{theorem}

{\it Proof.}
Choose $w\in \Z$. Let $((M,0),H,\nnn,P)$ be the $(trTLEP(w))$-structure
which corresponds to the Frobenius type structure on $K\to (M,0)$
by theorem \ref{t4.2}.
Its germ over $\cnnn$ is a $(TEP(w))$-structure with all the properties
in theorem \ref{t2.5}. Consider a universal unfolding of this germ
with base space $(\www M,0)= (M\times \C^l,0)$.
It extends uniquely to a $(trTLEP(w))$-structure 
$((\www M,0),\www H,\www\nnn,P)$ which unfolds $((M,0),H,\nnn,P)$.
This follows from the rigidity of logarithmic poles (lemma \ref{t2.6} (a)
and remark \ref{t3.3} (iii)) and from lemma \ref{t3.2}.

Let $((\www M,0),\www K, \www \nnn^r,\www C,\www \UU,\www \VV,\www g)$
be the Frobenius type structure which corresponds to this
$(trTLEP(w))$-structure by theorem \ref{t4.2}. There is a canonical 
isomorphism from the Frobenius type structure on $K$ to the restricted
one on $\www K|_{(M\times \{0\},0)}$. 
Let $\www \zeta\in \www K_0$ be the image of $\zeta\in K_0$.
It extends to a unique $\www \nnn^r$-flat section $\www v_1\in \OO(\www K)$.
The map
\begin{eqnarray}\label{4.12}
v: \TT_{\www M,0}\to \OO(\www K),\ X\mapsto -C_X\www v_1
\end{eqnarray}
is an isomorphism. It allows to shift the structure on $\www K$ to structure
on $T\www M$. Define
\begin{eqnarray}\label{4.13}
\nnn^v:=v^*\www \nnn^r,\ g^v :=v^*\www g,\ e:= v^{-1}(\www v_1),\  
E:= v^{-1}(\www \UU(\www v_1)).
\end{eqnarray}
The connection $\nnn^v$ on $T\www M$ is flat with $\nnn^v g^v=0$ and
$\nnn^v e=0$. Applying \eqref{t4.2} to $\www v_1$ shows that 
$\nnn^v$ is torsion free; so it is the Levi--Civita connection of $g^v$.
One defines a commutative and associative multiplication $\circ$ on 
$T\www M$ by
\begin{eqnarray}\label{4.14}
v(X\circ Y)= - C_X v(Y) = C_XC_Y \www v_1.
\end{eqnarray}
Then $e$ is the unit field. The potentiality condition follows from
\eqref{4.2}.
It rests to prove $\Lie_E(\circ)=\circ$ and $\Lie_E(g)=(2-d)\cdot g$.
We refer to \cite[Theorem 5.12]{He2} for details. The calculations
use \eqref{4.2}, \eqref{4.3}, \eqref{4.6}, $\www \nnn^r \www v_1=0$
and $\www\VV \www v_1=\frac{d}{2}\www v_1$.
They show especially
\begin{eqnarray}\label{4.15}
v(\nnn^v_X E) = (\www\VV + \frac{2-d}{2}\id )v(X) \mbox{ \ for }
X\in \TT_{M,0}.
\end{eqnarray}
One obtains a germ of a Frobenius manifold $((\www M,0),\circ, e,E,g^v)$.
Each step in its construction is essentially unique.
\hfill $\qed$

\begin{remarks}\label{t4.6}
(i) Theorem \ref{t4.5} is reduced with theorem \ref{t2.5} to the case 
when $M=\www M$ and when the map $C_\bullet \zeta:T_0M\to K_0$
is an isomorphism. This case was formulated
by Sabbah \cite[Theorem (4.3.6)]{Sab1} \cite[Th\'eor\`eme VII.3.6]{Sab2},
and independently by Barannikov \cite{Ba1}\cite{Ba2}\cite{Ba3}.
He called a major part of the initial data 
{\it semi-infinite variation of Hodge structures} (cf. remark \ref{3.3} (ii)).
Theorem \ref{t4.5} in the case $M=\www M$ is also implicit in the 
construction in singularity theory \cite{SK}\cite{SM}.

(ii) A Frobenius type structure with $M=\{0\}$ is simply a vector space
$K$ with a pairing $g$ and two endomorphisms $\UU$ and $\VV$ which
satisfy \eqref{4.5} and \eqref{4.6}. Then the condition (IC) in
theorem \eqref{4.5} is empty, (GC) must be satisfied by $\UU$ alone,
and (GC) and (EC) together are still more restrictive.
For example, this situation is satisfied at a point of a Frobenius manifold
where the multiplication with the Euler field is semisimple with 
different eigenvalues. This case was first considered by Dubrovin
\cite[Lecture 3]{Du}.

(iii) The case of a Frobenius type structure with $\UU=0$ can be considered
as opposite to the case in (ii). It will be discussed in chapter \ref{s5}.

(iv) One can define $(trTLP(w))$-structures and Frobenius type structures
without operators $\UU$ and $\VV$ \cite[5.2]{He2}. Omitting the corresponding
parts in theorem \ref{t4.2}, one obtains a one-to-one correspondence
between them. One can define Frobenius manifolds without Euler field.
The analogues of lemma \ref{t4.4} and theorem \ref{t4.5} hold.
This follows with the remarks \ref{t3.3}. But now the generation condition
(GC) requires $\dim M>0$.

(v) All the structures in chapters \ref{s2} to \ref{s4} were convergent
with respect to parameters $(t_1,...,t_m)\in (\C^m,0)\cong (M,0)$.
One can formulate everything in structures which are formal in these
parameters.

(vi) A. Kresch \cite[Theorem 1]{Kr} proved a strong reconstruction theorem
(existence and uniqueness)
for formal germs of Frobenius manifolds without Euler field and with
some additional properties typical for quantum cohomology.
It strengthens the first reconstruction theorem in \cite[Theorem 3.1]{KM},
which establishes only uniqueness.
It looks as if his result is the special case of the analogue of
theorem \ref{t4.5} without Euler field and with $(M,0)\subset (\www M,0)$
being the small quantum cohomology space. The collection
of $N(\beta,d)$ in \cite[Theorem 1]{Kr} gives the Higgs field $C$ on 
$T\www M|_M$,
the conditions on $A_1$ and $A$ give the generation condition (GC)'
and the property $C_XC_Y=C_YC_X$ for $X,Y\in \tm$ of the 
Higgs bundle. The other conditions on
Frobenius type structures (without $\UU$ and $\VV$) seem to be built-in.
\end{remarks}

\section[$H^2$-generated variations of filtrations]
{$H^2$-generated variations of filtrations and Frobenius manifolds}\label{s5}
\setcounter{equation}{0}

\noindent
In this chapter the special case of the construction theorem \ref{t4.5}
for Frobenius manifolds is studied when the Frobenius type structure
satisfies $\UU=0$. 
A Frobenius type structure with $\UU=0$ 
is equivalent to a variation of filtrations
with Griffiths transversality and additional structure (lemma \ref{t5.1}).
Typical cases are variations of polarized Hodge structures with
a condition, which is called $H^2$-generation condition, motivated by
quantum cohomology. The definitions \ref{t5.3} and \ref{t5.4} present
the relevant notions, theorem \ref{t4.5} reformulates the construction
theorem in the case $\UU=0$.

\begin{lemma}\label{t5.1}
(a) The structures in $(\alpha)$ and $(\beta)$ are equivalent.

$(\alpha)$ A Frobenius type structure $((M,0),K,\nnn^r,C,\UU,\VV,g)$ together
with an integer $w$ such that $\UU=0$ and $\VV$ is semisimple with eigenvalues
in $\frac{w}{2}+\Z$.

$(\beta)$ A tuple $((M,0),K,\nnn,F^\bullet,U_\bullet,w,S)$. Here
$K\to (M,0)$ is a germ of a holomorphic vector bundle; $\nnn$ is a flat
connection on $K$; $F^\bullet$ is a decreasing filtration by 
germs of holomorphic subbundles $F^p\subset K$, $p\in \Z$, which
satisfies Griffiths transversality
\begin{eqnarray}\label{5.1}
\nnn: \OO(F^p)\to \Omega_{M,0}^1\otimes \OO(F^{p-1});
\end{eqnarray}
$U_\bullet$ is an increasing filtration by flat subbundles $U_p\subset K$
such that
\begin{eqnarray}\label{5.2}
H=F^p\oplus U_{p-1} = \bigoplus_q F^q\cap U_q;
\end{eqnarray}
$w\in \Z$, and $S$ is a $\nnn$-flat, $(-1)^w$-symmetric, nondegenerate
pairing on $K$ with
\begin{eqnarray}\label{5.3}
S(F^p,F^{w+1-p})=0,\\
S(U_p,U_{w-1-p})=0.\label{5.4}
\end{eqnarray}

(b) One passes from $(\alpha)$ to $(\beta)$ by defining
\begin{eqnarray}\label{5.5}
\nnn &:=& \nnn^r + C,\\
F^p &:=& \bigoplus_{q\geq p} \ker (\VV-(q-\frac{w}{2})\id:K\to K), 
\label{5.6}\\
U_p &:=& \bigoplus_{q\leq p} \ker (\VV-(q-\frac{w}{2})\id:K\to K), 
\label{5.7}\\
S(a,b) &:=& (-1)^p g(a,b) \mbox{ for } a\in\OO(F^p\cap U_p), b\in \OO(K).
\label{5.8}
\end{eqnarray}
\end{lemma}

{\it Proof.}
First we prove part (b). The connection $\nnn^r$ and the 
Higgs field $C$ are maps 
$\OO(K)\to \Omega^1_{M,0}\otimes \OO(K)$. 
They can be extended canonically to
maps $\Omega^1_{M,0}\otimes \OO(K) \to \Omega^2_{M,0}\otimes \OO(K) $.
Then the flatness of $\nnn^r$ means $(\nnn^r)^2 =0$, the Higgs field
$C$ satisfies $C^2=0$, and the potentiality condition \eqref{4.2}
means $\nnn^r(C) := \nnn^r\circ C+C\circ \nnn^r = 0$. Therefore
$\nnn^2= (\nnn^r+C)^2=0$, the connection $\nnn$ is flat.
The filtrations $F^\bullet $ and $U_\bullet$ obviously satisfy \eqref{5.2}
and
\begin{eqnarray}\label{5.9}
F^p\cap U_p =\ker (\VV-(p-\frac{w}{2})\id: K\to K).
\end{eqnarray}
The connection $\nnn^r$ maps $\OO(F^p\cap U_p)$ to itself because
$\VV$ is $\nnn^r$-flat. Because of $\UU=0$ the condition \eqref{4.3}
is $[C,\VV]=C$. Equivalent is that $C_X$, $X\in \TT_{M,0}$, maps
$\OO(F^p\cap U_P)$ to $\OO(F^{p-1}\cap U_{p-1})$.
Therefore $U_\bullet$ is $\nnn$-flat and $F^\bullet$ satisfies
Griffiths transversality.

The condition \eqref{4.6} says that for $a\in \OO(F^p\cap U_p)$,
$b\in \OO(F^q\cap U_q)$
\begin{eqnarray}\label{5.10}
S(a,b) = (-1)^pg(a,b) = 0 \mbox{ if } p+q\neq w.
\end{eqnarray}
Therefore $S$ is $(-1)^w$-symmetric and satisfies \eqref{5.3} and 
\eqref{5.4}. It is $\nnn$-flat because for $X\in \TT_{M,0}$ and
$\nnn^r$-flat sections $a\in \OO(F^p\cap U_p)$, $b\in \OO(K)$
\begin{eqnarray} \nonumber
(\nnn_XS)(a,b) &=& X\, S(a,b)-S(\nnn^r_Xa+C_X a,b) 
- S(a,\nnn^r_Xb+C_X b) \\
&=& 0- (-1)^{p+1}g(C_X a,b) - (-1)^pg(a,C_X b)=0.
\label{5.11}
\end{eqnarray}
This shows part (b). One passes from $(\beta)$ to $(\alpha)$ as follows.
One defines the endomorphism $\VV$ by \eqref{5.9}, the pairing $g$
by \eqref{5.8}, and one decomposes $\nnn$ into $\nnn^r$ and $C$ such 
that $\nnn^r$ maps $\OO(F^p\cap U_p)$ to itself and $C_X$ for $X\in \TT_{M,0}$
maps $\OO(F^p\cap U_p)$ to $\OO(F^{p-1}\cap U_{p-1})$.
One easily checks all conditions of a Frobenius type structure.
\hfill $\qed$

\begin{remarks}\label{t5.2}
(i) If one adds in lemma \ref{t5.1} (a) $(\beta)$ real structure with 
suitable conditions then one obtains a germ of a variation of 
polarized Hodge structures of weight $w$.

(ii) A Frobenius type structure is equivalent to a $(trTLEP(w))$-structure
(theorem \ref{t4.2}), which is composed of a $(TEP(w))$-structure
at $\cnnn$ and an $(LEP(w))$-structure at $(\P^1,\infty)\times (M,0)$.
One can refine lemma \ref{t5.1} (a). There is a correspondence between
$(TEP(w))$-structures with $\UU=0$ and monodromy $(-1)^w\id$ on the 
one hand and germs of variations of filtrations with Griffiths 
transversality and a pairing on the other hand
\cite[Corollary 7.14]{He2}.
Under this correspondence an $(LEP(w))$-structure corresponds to 
a trivial, flat variation of filtrations. Putting these together,
a $(trTLEP(w))$-structure corresponds to the structure in lemma \ref{5.1} (a)
$(\beta)$.

But the variation of filtrations $F^\bullet$ which corresponds to a 
$(TEP(w))$-structure $(H\to \cnnn,\nnn)$ does not live on the bundle
$K=H|_{\nmmn}$. It lives on a flat bundle on $(M,0)$ whose fibers
are all isomorphic to the space of global flat manyvalued sections
in $H|_\csnnn$. The bundle $K$ is only isomorphic to 
$\bigoplus_p F^p/F^{p+1}$. Only in the case of a second filtration
$U_\bullet$ and a splitting \eqref{5.2} one obtains an isomorphism
between $K$ and the flat bundle on which $F^\bullet$ lives.
\end{remarks}

\begin{definition}\label{t5.3}
(a) A germ of an {\it $H^2$-generated variation of filtrations of weight $w$}
is a tuple $((M,0),K,\nnn,F^\bullet,w)$ with the following properties.
$w$ is an integer, $K\to (M,0)$ is a germ of a holomorphic vector bundle
with flat connection $\nnn$ and variation of filtrations $F^\bullet$
which satisfies Griffiths transversality \eqref{5.1}
and
\begin{eqnarray}\label{5.12}
&& 0=F^w\subset F^{w-1}\subset ... \subset  K \\
&& \rk F^{w-1}=1,\ \ \rk F^{w-2}=1+\dim M \geq 2.\label{5.13}
\end{eqnarray}
Griffiths transversality and flatness of $\nnn$ give a Higgs field $C$
on the graded bundle $\bigoplus_p F^p/F^{p+1}$ with commuting 
endomorphisms 
\begin{eqnarray}\label{5.14}
C_X = [\nnn_X]:\OO(F^p/F^{p+1}) \to \OO(F^{p-1}/F^p) \mbox{ for }
X\in \TT_{M,0}.
\end{eqnarray}
{\it $H^2$-generation condition:}
the whole module $\bigoplus_p \OO(F^p/F^{p+1})$ is generated by
$\OO(F^{w-1})$ and its images under iterations of the maps
$C_X$, $X\in \TT_{M,0}$.

(b) A pairing and an opposite filtration for an $H^2$-generated 
variation of filtrations $((M,0),K,\nnn,F^\bullet)$ of weight $w$ are
a pairing $S$ and a filtration $U_\bullet$ as in lemma \ref{t5.1} (a)
$(\beta)$.
\end{definition}

\begin{definition}\label{t5.4}
An $H^2$-generated germ of a Frobenius manifold of weight $w\in \N_{\geq 3}$
is a germ $((M,0),\circ,e,E,g)$ of a Frobenius manifold with the 
properties (I) and (II) below and with
\begin{eqnarray}\label{5.15}
E|_{t=0}=0.
\end{eqnarray}
Let $\nnn^g$ be the Levi--Civita connection of $g$. The endomorphism
$\nnn^gE:\TT_{M,0}\to \TT_{M,0}$, $X\mapsto \nnn^g_XE$, acts
on the space of $\nnn^g$-flat vector fields. In particular
$e\in \ker(\nnn^gE-\id)$.

(I) It acts semisimply with eigenvalues $\{1,0,...,-(w-3)\}$.

It turns out that then the multiplication on the algebra $T_0M$ respects
the grading 
\begin{eqnarray}\label{5.16}
T_0M = \bigoplus_{p=0}^{w-2} 
\ker(\nnn^g E - (1-p)\id:T_0M\to T_0M)
\end{eqnarray}
and that $\Lie_E(g)=(4-w)\cdot g$.

(II) $H^2$-generation condition: The algebra $T_0M$ is generated by
$\ker (\nnn^gE:T_0M\to T_0M)$.
\end{definition}

\begin{remarks}\label{t5.4b}
(i) Properties (I) and \eqref{5.15} hold for even-dimensional 
quantum cohomology of Calabi--Yau manifolds of complex dimension $w-2$.
This follows from the vanishing of the canonical class and the
standard explicit formulas for the Euler field in quantum cohomology.
Generally, our terminology involving ``$H^2$-generation'' is motivated
by quantum cohomology, for which $\ker \nnn^gE =H^2$.
For a more extended discussion of special properties of quantum cohomology
Frobenius manifolds, see 
\cite[1.3 and 1.4, in particular Definition 1.4.1]{Man1}.

(ii) The uniqueness statement in theorem \ref{t5.5} that a structure as
in ($\gamma$) is determined by a structure as in ($\beta$) is essentially
a special case of that in \cite[Theorem 3.1 and 3.1.1 a) and b)]{KM}. 
A refined version of it (cf. lemma \ref{t8.2})
was used already in the proof of \cite[Theorem 6.5]{Ba1}.
But the existence statement that any structure as in ($\beta$) gives
rise to a structure as in ($\gamma$) is new.
\end{remarks}

\begin{theorem}\label{t5.5}
There is a one-to-one correspondence between the structures in 
$(\alpha)$, $(\beta)$ and $(\gamma)$.

$(\alpha)$ A Frobenius type structure $((M,0),K,\nnn^r,C,\UU,\VV,g)$
with $\UU=0$ and with a fixed vector $\zeta\in K_0$ which satisfies
the conditions (IC), (GC) and (EC) in theorem \ref{t4.5}.

$(\beta)$ A germ of an $H^2$-generated variation of filtrations
$((M,0),K,\nnn,F^\bullet,w,S, U_\bullet)$ of weight $w\in \N_{\geq 3}$
with pairing and opposite filtration and with a fixed generator
$\zeta\in (F^{w-1})_0\subset K_0$.

$(\gamma)$ An $H^2$-generated germ of a Frobenius manifold
$((\www M,0),\circ,e,E,\www g)$ of weight $w\in \N_{\geq 3}$.

One passes from $(\alpha)$ to $(\beta)$ by lemma \ref{t5.1} (b) and from 
$(\alpha)$ to $(\gamma)$ by theorem \ref{t4.5}.
One passes from $(\gamma)$ to $(\alpha)$ by defining 
\begin{eqnarray}\label{5.17}
M:= \{t\in \www M\ |\ E|_t=0\},
\end{eqnarray}
$K:=T\www M|_{(M,0)}$ with the canonical Frobenius type structure,
and $\zeta:=e|_0$. The eigenvector condition (EC) is 
$\VV\zeta = \frac{w-2}{2}\zeta$.
\end{theorem}

{\it Proof.}
Let us start with $(\alpha)$. We have to show that for some $w\in \N_{\geq 3}$
the endomorphism $\VV+\frac{w}{2}\id$ is semisimple with eigenvalues
in $\{1,2,...,w-1\}$ and with $\VV \zeta = \frac{w-2}{2}\zeta$.

The eigenvector condition (EC) says $\VV\zeta = \frac{d}{2}\zeta$
for some $d\in \C$. Condition \eqref{4.3} reads as $[C,\VV]=C$.
This together with the generation condition (GC) shows that $\VV$
acts semisimply on $K_0$ with eigenvalues
in $\frac{d}{2}+\Z_{\leq 0}$ and that $\ker(\VV-\frac{d}{2}\id)=\C\cdot\zeta$.
The injectivity condition (IC) implies $\dim \ker(\VV-(\frac{d}{2}-1)\id)
=\dim M>0$.
Condition \eqref{4.6} tells that the eigenvalues of $\VV$ are in 
$(\frac{d}{2}+\Z_{\leq 0})\cap -(\frac{d}{2}+\Z_{\leq 0})$.
Therefore $d\in \N$. Define $w:=d+2$.

The conditions (IC) and (GC) show that one passes from $(\alpha)$
to $(\beta)$ by lemma \ref{t5.1} (b). It is also clear that one
can pass back.

Theorem \ref{t4.5} applied to $(\alpha)$ gives a Frobenius manifold
$((\www M,0),\circ,e,E,\www g)$, an embedding $i:(M,0)\to (\www M,0)$
and an isomorphism $j:K\to T\www M|_{i(M)}$ of Frobenius type structures.
It maps $\VV$ to the restriction of $\nnn^gE + \frac{2-d}{2}\id$ on 
$T\www M|_{i(M)}$ \eqref{4.15}. Therefore one obtains an 
$H^2$-generated germ of a Frobenius manifold of weight $w$.
In suitable flat coordinates the Euler field of this Frobenius manifold
takes the form
\begin{eqnarray}\label{5.18}
E = \sum_{i=1}^{\dim \www M} d_it_i\frac{\paa}{\paa t_i}
\end{eqnarray}
with $d_i\in \{1,0,...,-(w-3)\}$ and 
\begin{eqnarray}\label{5.19}
\sharp (i\ |\ d_i=0)=\dim \ker(\nnn^gE) = \dim \ker(\VV-(\frac{d}{2}-1)\id)
=\dim M.
\end{eqnarray}
Therefore $i(M)= \{t\in \www M\ |\ E|_t=0\}$.
One sees also easily that one passes from $(\gamma)$ to $(\alpha)$ as
described in the theorem.
\hfill $\qed$

\bigskip
A distinguished class of $H^2$-generated variations of filtrations are
variations of Hodge filtrations associated to certain families of 
homogeneous polynomials. We will discuss them and the corresponding
$H^2$-generated Frobenius manifolds in the chapters \ref{s6}, \ref{s7}
and \ref{s8}. The following example shows that it is easy to
construct abstract $H^2$-generated variations of filtrations.
Therefore one has a lot of freedom in constructing $H^2$-generated
Frobenius manifolds.

\begin{example}\label{t5.6}
Consider $(M,0):=(\C,0)$ with coordinate $t$, the trivial bundle
$H:= \C^{w-1}\times (M,0)\to (M,0)$ for some $w\in \N_{\geq 3}$ 
with standard basis $v_1,...,v_{w-1}$ of sections, a pairing $S$ 
with $S(v_p,v_q):=(-1)^p\delta_{p,w-q}$, and filtrations $F^\bullet$ and
$U_\bullet$ with
\begin{eqnarray}\label{5.20}
\OO(F^p):= \bigoplus_{q\geq p}\OO_{\C,0}\cdot v_{w-q},\ \ 
\OO(U_p):=\bigoplus_{q\leq p}\OO_{\C,0}\cdot v_{w-q}.
\end{eqnarray}
Choose {\it any} invertible functions 
$b_2,...,b_{\left[ \frac{w-1}{2}\right]}\in \OO^*_{\C,0}$, define 
\begin{eqnarray}\label{5.21}
b_1:=1, \ b_{w-1}:=0, \ b_k:=b_{w-1-k}\mbox{ for }
k=\left[ \frac{w-1}{2}\right]+1,...,w-2
\end{eqnarray}
and define a 
connection $\nnn$ on $H$ by 
\begin{eqnarray}\label{5.22}
\nnn_{\paa/\paa t} v_i:= b_iv_{i+1}.
\end{eqnarray}
Then $((M,0),K,\nnn,F^\bullet,w,S,U_\bullet)$ is a germ of an $H^2$-generated
variation of filtrations of weight $w$ with pairing  and opposite filtration.
Moreover, one can see that a second tuple
$(b'_2,...,b'_{\left[ \frac{w-1}{2}\right]})$ of functions yields isomorphic
data only if $b_i$ and $b'_i$ coincide up to multiplication by a constant.
The big freedom in constructing $H^2$-generated Frobenius manifolds of 
dimension $w-1$ 
is in striking difference to the semisimple case, where one has only
finitely many parameters \cite[II 3.4.3]{Man2}.
\end{example}

\section{Hypersurfaces in $\P^n$}\label{s6}
\setcounter{equation}{0}

\noindent
In the case of certain families of smooth hypersurfaces in $\P^n$ one
obtains variations of Hodge structures which are the prototype of 
$H^2$-generated variations of filtrations (definition \ref{t5.3}). 
This is a simple consequence
of Griffiths' description of the Hodge filtration on the primitive
part of the middle cohomology of a smooth hypersurface in $\P^n$ in terms
of rational differential forms on $\P^n$ \cite{Gr}.

Fix a degree $d\in \N$ and denote by $\C[x]^{(q)}$ for 
$q\in \frac{1}{d}\Z_{\geq 0}$
the space of homogeneous polynomials in $\C[x_0,...,x_n]=\C[x]$ of
degree $d\cdot q$. 
Consider a polynomial $f\in \C[x]^{(1)}$ with isolated singularity at 0.
The grading on $\C[x]$ induces a grading on the Jacobi algebra
\begin{eqnarray}\label{6.1}
R_f:= \C[x]/\left(\frac{\paa f}{\paa x_0},...,\frac{\paa f}{\paa x_n}\right)
\end{eqnarray}
of $f$ with subspaces $R_f^{(q)}$, $q\in \frac{1}{d}\Z_{\geq 0}$. 
The primitive part of the middle cohomology $H^{n-1}(X,\C)$
of the smooth hypersurface $X:=\oooo{f^{-1}(0)}\subset \P^n$ is denoted by
$H^{n-1}_{prim}(X)$, its Hodge filtration by $F^\bullet_{prim}\subset 
H^{n-1}_{prim}(X)$. The primitive cohomology also comes equipped 
with a polarizing form $S$.

Now consider a family of polynomials $F_t\in \C[x]^{(1)}$, $t\in M_0$,
with isolated singularities at 0, where $M_0$ is a smooth parameter space
with coordinates $(t_\alpha)$ such that for each $t\in M_0$ the map
\begin{eqnarray}\label{6.2}
\aaa :T_tM_0 \to R_{F_t}^{(1)},\ \ \frac{\paa}{\paa t_\alpha} \mapsto
\left[ \frac{\paa F_t}{\paa t_\alpha}\right] 
\end{eqnarray}
is an isomorphism.
The bundle $H:=\bigcup_{t\in M_0} H^{n-1}_{prim}(X_t)$ comes equipped with 
a real subbundle, a flat connection $\nnn$, the flat pairing $S$ and 
a variation of Hodge filtrations $F^\bullet_{prim}$. Together they form
a variation of polarized Hodge structures of weight $n-1$.

\begin{theorem}\label{t6.1}
If $\frac{n+1}{d}\in\N$ then the tuple $(M_0,H,\nnn,S,\www F^\bullet :=
F_{prim}^{\bullet-2+(n+1)/d})$ is an $H^2$-generated variation of filtrations
of weight $n+3-2(n+1)/d$ with pairing.
\end{theorem}

{\it Proof.}
The space $\Omega_{alg}^{n+1}=\C[x]\ddd x_0...\ddd x_n$ of algebraic
differential forms on $\C^{n+1}$ is graded with subspaces 
$(\Omega_{alg}^{n+1})^{(q)} = \C[x]^{(q-(n+1)/d)}\ddd x_0...\ddd x_n$ for
$q\in \frac{1}{d}\Z_{\geq n+1}$. For $f\in \C[x]^{(1)}$ with 
isolated singularity the quotient
\begin{eqnarray}\label{6.3}
\Omega_f := \Omega_{alg}^{n+1}/\ddd f\land \Omega_{alg}^n
\end{eqnarray}
carries an induced grading with subspaces $\Omega_f^{(q)}$.
It is a graded module of the graded algebra $R_{f}$, and it is 
a free module of rank 1 of $R_f$.

Following Griffiths \cite{Gr}, one obtains for $q\in \Z_{\geq (n+1)/d}$
a canonical isomorphism
\begin{eqnarray}\label{6.4}
\rrr_q:\Omega_f^{(q)} \to F^{n-q}_{prim}/ F^{n-q+1}_{prim}
\end{eqnarray}
in the following way.
Let $\EE := \sum_{i=0}^n x_i\frac{\paa}{\paa x_i}$ be the Euler field
on $\C^{n+1}$. For $\omega\in (\Omega_{alg}^{n+1})^{(q)}$, 
$q\in \Z_{\geq (n+1)/d}$, consider the form $i_\EE(\frac{\omega}{f^q})$,
which one obtains from $\frac{\omega}{f^q}$ by contraction with the
Euler field. It
extends to a rational form on $\P^{n+1}=\C^{n+1}\cup \P^n$ with a pole
of order $\leq q$ along ${f^{-1}(0)}\cup X$.
We denote the restriction of this form to $\P^n$ by 
$\oooo{i_\EE (\frac{\omega}{f^q})}$.
It induces a class $\left[ \oooo{i_\EE(\frac{\omega}{f^q})}\right]
\in H^n(\P^n-X,\C)$. There is a residue map
\begin{eqnarray}\label{6.5}
\Res : H^n(\P^n-X,\C)\to H^{n-1}(X,\C),
\end{eqnarray}
dual to a tube map in the homologies. One defines
\begin{eqnarray}\label{6.6}
\rho_q : (\Omega_{alg}^{n+1})^{(q)} \to H^{n-1}(X,\C),\\
   \omega\mapsto \Res\left[ \oooo{i_\EE\left(\frac{\omega}{f^q}\right)}\right]
\cdot \frac{(q-1)!}{d}.\nonumber 
\end{eqnarray}
By \cite{Gr}, the image of $\rho_q$ is $F^{n-q}_{prim}\subset H^{n-1}_{prim}(X)$
and the preimage of $F_{prim}^{n-q+1}$ 
is $(\ddd f\land \Omega_{alg}^n)^{(q)}$. This gives the isomorphism $\rrr_q$.

Now for a family of polynomials $F_t\in \C[x]^{(1)}$, $t\in M_0$, with
isolated singularities at 0 the infinitesimal variation of Hodge structures
can be calculated. For each $t\in M_0$ it is a set of commuting linear maps
\begin{eqnarray}\label{6.7}
C_{\paa/\paa  t_\alpha} = [\nnn_{\paa/\paa  t_\alpha}]:
F^p_{prim}/F^{p+1}_{prim} \to F_{prim}^{p-1}/F_{prim}^{p}.
\end{eqnarray}
The calculation (cf. \cite[2.2]{Do})
\begin{eqnarray}\label{6.8}
\frac{\paa}{\paa t_\alpha} \frac{\omega}{F_t^q} = 
(-q)\frac{\frac{\paa F_t}{\paa t_\alpha} \cdot \omega}{F_t^{q+1}}
\end{eqnarray}
shows for $[\omega]\in \Omega_f^{(q)}$
\begin{eqnarray}\label{6.9}
C_{\paa/\paa t_\alpha}  \rrr_q([\omega]) = \rrr_{q+1}\left( -\aaa 
(\frac{\paa}{\paa t_\alpha})\cdot [\omega]\right).
\end{eqnarray}
Now one observes two facts: (i) The subring
$\bigoplus_{q\in \Z_{\geq 0}}R_{F_t}^{(q)}$ of $R_{F_t}$ is multiplicatively
generated by $R_{F_t}^{(1)}$, because any monomial in $\C[x]^{(q)}$ for
$q\in \Z_{\geq 1}$ is a product of monomials in $\C[x]^{(1)}$.\\
(ii) If $\frac{n+1}{d}\in \N$ then $\bigoplus_{q\in \Z}\Omega_{F_t}^{(q)}$
is a free module of rank 1 of the ring
$\bigoplus_{q\in \Z_{\geq 0}}R_{F_t}^{(q)}$.

This shows that for a family of polynomials with \eqref{6.2} and 
$\frac{n+1}{d}\in \N$ the $H^2$-generation condition in definition \ref{t5.3}
is satisfied.
\hfill $\qed$

\bigskip
Via theorem \ref{t5.5} the variations of Hodge filtrations in 
theorem \ref{t6.1} together with chosen opposite filtrations
give rise to Frobenius manifolds. All of them can be
identified with submanifolds of Frobenius manifolds which arise in 
singularity theory, see chapter \ref{s7}.
Those with $d=n+1$ can be identified with submanifolds of Frobenius (super)
manifolds in the Barannikov--Kontsevich construction, see chapter \ref{s8}.
The hypersurfaces $X=\oooo{f^{-1}(0)}\subset \P^n$ are Calabi--Yau if and
only if $d=n+1$.

Much of the preceding discussion generalizes to the case of quasihomogeneous
singularities, but not all. Lemma \ref{t6.2} gives an example where the
$H^2$-generation condition fails to hold. 

A weight system $(w_0,...,w_n)$ with $w_i\in \Q\cap (0,\frac{1}{2}]$ 
induces a new grading on $\C[x_0,...,x_n]$ whose subspaces we also 
denote by $\C[x]^{(q)}$. A monomial $x_0^{i_0}...x_n^{i_n}$ is in 
$\C[x]^{(q)}$ for $q\in \Q_{\geq 0}$ if $\sum i_jw_j=q$. 
The Jacobi algebra $R_f$ of $f\in \C[x]^{(1)}$ is defined as above and 
inherits a grading with subspaces $R_f^{(q)}$.

\begin{lemma}\label{t6.2}
Consider $n=5$ and $(w_0,...,w_5)=\frac{1}{9}(1,1,1,2,2,2)$ and any weighted
homogeneous polynomial $f\in \C[x]^{(1)}$ with isolated singularity at 0.
The subspace of $R_f^{(2)}$ which is generated by the set
$R_f^{(1)}\cdot R_f^{(1)}$ has codimension 1 in $R_f^{(2)}$.
\end{lemma}

{\it Proof.}
A monomial in $\C[x]^{(2)}$ is a product of monomials in $\C[x]^{(1)}$ 
if and only if it contains $x_0,x_1$ or $x_2$. Therefore it is sufficient
to show that the space
\begin{eqnarray}\label{6.10}
\C[x_3,x_4,x_5]^{(2)}\cap \left(\frac{\paa f}{\paa x_0},...,
\frac{\paa f}{\paa x_n}\right)
\end{eqnarray}
has codimension 1 in $\C[x_3,x_4,x_5]^{(2)}$.
Notice $\frac{\paa f}{\paa x_i}|_{\{x_0=x_1=x_2=0\}} =0$ for $i=3,4,5$
and define $\www f_i (x_3,x_4,x_5):= 
\frac{\paa f}{\paa x_i}|_{\{x_0=x_1=x_2=0\}} $ 
for $i=0,1,2$. The ideal $(\www f_0,\www f_1,\www f_2)\subset \C[x_3,x_4,x_5]$
has an isolated zero at 0. The dimension of 
$\left(\C[x_3,x_4,x_5]/(\www f_0,\www f_1,\www f_2)\right)^{(2)}$
is independent of the choice of $\www f_0,\www f_1,\www f_2$
as long as the ideal $(\www f_0,\www f_1,\www f_2)$ has an isolated zero
at 0. The choice $\www f_0=x_3^4$, $\www f_1=x_4^4$, $\www f_2=x_5^4$
shows that the dimension is 1.
\hfill $\qed$

\section{Frobenius manifolds for hypersurface singularities}\label{s7}
\setcounter{equation}{0}

\noindent
In theorem \ref{t7.3} the Frobenius manifolds which one obtains from
theorem \ref{6.1} combined with theorem \ref{5.5} will be identified with 
submanifolds of Frobenius manifolds in singularity theory.
For each holomorphic function germ $f:(\C^{n+1},0)\to (\C,0)$ 
with an isolated singularity at 0 the base space $M_\mu\subset \C^\mu$
of a semiuniversal unfolding can be equipped with the structure
of a Frobenius manifold \cite{SK}\cite{SM}. A detailed account
is given in \cite{He1}. 
In \cite[ch. 8]{He2} the construction is recasted as the construction
of a $(trTLEP(w))$-structure on $M_\mu$ and subsequent application of
a special case of theorem \ref{t4.5}.
We restrict now to quasihomogeneous singularities
and recall some facts from \cite{He1}.

Let $(w_0,...,w_n)$ be a weight system with $w_i\in \Q\cap (0,\frac{1}{2}]$
and $f\in \C[x]^{(1)}$ (notation from the end of chapter 6) 
a weighted homogeneous polynomial with isolated singularity at 0 and 
Milnor number $\mu$. 
From the grading of the Jacobi algebra $R_f$ one obtains the exponents
$\alpha_1,...,\alpha_\mu$, rational numbers with
\begin{eqnarray}\label{7.1}
\alpha_1\leq ... \leq \alpha_\mu,\ \ \alpha_1=\sum w_i,\ \  
\sharp (i\ |\ \alpha_i=\alpha) =\dim R_f^{(\alpha-\alpha_1)}.
\end{eqnarray}
Choose polynomials $m_i\in \C[x]^{(\alpha_i-\alpha_1)}$ which represent a
basis of the Jacobi algebra and such that $m_1=1$. The function
\begin{eqnarray}\label{7.2}
F(x_0,...,x_n,t_1,...,t_\mu) = f+\sum_{i=1}^\mu t_im_i
\end{eqnarray}
is a semiuniversal unfolding of $f$. 
It should be seen as a family of functions $F_t$ with 
parameter $t=(t_1,...,t_\mu)\in M_\mu
\subset \C^\mu$. Here $M_\mu$ is a suitable open neighborhood of 0.
The manifold $M_\mu$ comes equipped with the unit field 
$e=\frac{\paa }{\paa t_1}$, the Euler field
$E=\sum_{i=1}^\mu (1+\alpha_1-\alpha_i)t_i\frac{\paa}{\paa t_i}$, 
and a multiplication $\circ$ on the holomorphic tangent bundle.
The multiplication is induced from a canonical isomorphism of $T_tM_\mu$
with the direct sum of Jacobi algebras of the singularities of the
function $F_t$ for $t\in M_\mu$. 
The tuple $(M_\mu,\circ,e,E)$ is an F-manifold \cite[I \S5]{Man2}\cite[I]{He1}.

With the Gau{\ss}--Manin connection one can construct a metric $g$ on 
$M_\mu$ such that $(M_\mu,\circ,e,E,g)$ is a Frobenius manifold 
\cite{SK}\cite{SM}. In general the metric depends on a choice.
For the construction we refer to \cite{He1}. Here we merely explain
the choice and state the result.

Let $H^\infty$ be the space of global flat multivalued sections of the flat
cohomology bundle $\bigcup_{z\in \C^*}H^n(f^{-1}(z),\C)$.
It comes equipped with a real subbundle $H^\infty_\R$, a semisimple 
monodromy operator $h:H^\infty\to H^\infty$, a monodromy invariant
Hodge filtration $F^\bullet$ and a monodromy invariant pairing $S^\infty$
\cite[ch. 10]{He1}.
Define $H^\infty_\lambda := \ker(h-\lambda\id)\subset H^\infty$
and $H^\infty_{\neq 1} :=\bigcup_{\lambda\neq 1}H^\infty_\lambda $.

Then $(H^\infty_{\neq 1},H^\infty_\R\cap H^\infty_{\neq 1}, 
F^\bullet, S^\infty)$ and 
$(H^\infty_{1},H^\infty_\R\cap H^\infty_{1}, 
F^\bullet, S^\infty)$ are polarized Hodge structures of weight $n$ and $n+1$.
An increasing monodromy invariant filtration $U_\bullet$ on $H^\infty$
is called {\it opposite to } $F^\bullet$ if 
\begin{eqnarray}\label{7.3}
&& H^\infty = \bigoplus_pF^p\cap U_p,\\
&& S^\infty(H^\infty_{\neq 1}\cap F^p\cap U_p,
    H^\infty_{\neq 1}\cap F^q\cap U_q)=0 \mbox{ for }p+q\neq n,\label{7.4}\\
&& S^\infty(H^\infty_{1}\cap F^p\cap U_p,
    H^\infty_{1}\cap F^q\cap U_q)=0 \mbox{ for }p+q\neq n+1.\label{7.5}
\end{eqnarray}

\begin{theorem}\label{t7.1}\cite{SM}\cite[Theorem 11.1]{He1}
Any choice of an opposite filtration $U_\bullet$ induces an up to a scalar
unique metric $g$ on $M_\mu$ such that $(M_\mu,\circ,e,E,g)$ is a Frobenius 
manifold. The opposite filtration is uniquely determined by the metric.
\end{theorem}

Consider the submanifold
\begin{eqnarray}\label{7.6}
M:= \{t\in M_\mu\ |\ t_i=0 \mbox{ if }\alpha_i-\alpha_1\notin\Z\}
\subset M_\mu.
\end{eqnarray}
It does not depend on the choice of the coordinates $t_i$; 
that means, any choice
with $Et_i=(1+\alpha_1-\alpha_i)t_i$ gives the same submanifold $M$.
It parametrizes the semiquasihomogeneous deformations of $f$
by polynomials of integer degree.
Unit field $e$ and Euler field $E$ are tangent to $M$.
One can show that the multiplication $\circ$ on $TM_\mu$ restricts to
a multiplication on $TM$. If $\alpha_1=\sum w_i \in \frac{1}{2}\N$
then a much stronger result holds.

\begin{theorem}\label{t7.2}\cite[III 8.7.1]{Man2}
Suppose that $\alpha_1=\sum w_i \in \frac{1}{2}\N$.
For any metric $g$ as in theorem \ref{t7.1}, the submanifold $M$ with 
induced multiplication, metric, Euler field $E$ and unit field $e$ 
is a Frobenius manifold.
\end{theorem}

Going through the construction in \cite[11.1]{He1} one can even see that
the metric on $M$ depends only on $U_\bullet\cap H^\infty_{(-1)^{\alpha_1}}$
and that it determines $U_\bullet\cap H^\infty_{(-1)^{\alpha_1}}$ uniquely.

Finally we restrict to the case $(w_0,...,w_n)=\frac{1}{d}(1,...,1)$ for 
some $d\in \N$ with $\alpha_1=\frac{n+1}{d}\in \N$.
We use the notations in chapter \ref{s5}.
The manifold 
\begin{eqnarray}\label{7.7}
M_0:=\{t\in M_\mu\ |\ t_i =0 \mbox{ if }\alpha_i-\alpha_1\neq 1\}
\subset M
\end{eqnarray}
parametrizes the homogeneous polynomials $F_t$ in the unfolding
$F$. The map $\aaa$ in \eqref{6.2} is an isomorphism.
Theorem \ref{t6.1} applies to the family of functions $F_t$, $t\in M_0$.
As in chapter \ref{s6}, $X_t:= \oooo{F_t^{-1}(0)}\subset\P^n$
is the hypersurface in $\P^n$ defined by $F_t$.

\begin{theorem}\label{t7.3}
(a) There is a canonical isomorphism
\begin{eqnarray}\label{7.8}
\psi: \bigcup_{t\in M_0}H^\infty_1(F_t) \to 
\bigcup_{t\in M_0} H^{n-1}_{prim}(X_t)
\end{eqnarray}
of flat bundles with pairings $S^\infty$ and $S$. For each $t\in M_0$
it is a $(-1,-1)$ morphism of Hodge structures.

(b) The Frobenius manifold structures on $M$ from theorem \ref{t7.2}
and the Frobenius manifolds constructed in theorem \ref{t5.5} from the 
variation of filtrations $((M_0,0),H,\nnn,S, \www F^\bullet , \www U_\bullet)$
in theorem \ref{t6.1} with opposite filtration $\www U_\bullet$ 
on $H=\bigcup_{t\in M_0}H^{n-1}_{prim}(X_t)$ are pairwise isomorphic.
Two are isomorphic up to multiplication of the metric by a scalar
if and only if the opposite filtrations $U_\bullet\cap H^\infty_1$ 
on $H^\infty_1$ and $\www U_\bullet$ on $H$ satisfy
\begin{eqnarray}\label{7.9}
\psi( U_{\bullet -1+(n+1)/d}\cap H^\infty_1)= \www U_\bullet.
\end{eqnarray}
\end{theorem}

{\it Proof.}
(a) This is essentially well known. The following explanations
may be helpful.
Consider for some $z\in \C^*$ a fiber $F_t^{-1}(z)\subset \C^{n+1}$.
The hypersurface $X_t\subset\P^n$ is the part in $\P^n$
of the closure of $F_t^{-1}(z)\subset \C^{n+1}$ in $\P^{n+1}= \C^{n+1}\cup
\P^n$. Consider a tubular neighborhood $T(F_t^{-1}(z))$ of 
$X_t$ in $F_t^{-1}(z)$. There are canonical isomorphisms
\begin{eqnarray}\label{7.10}
H^\infty_1\leftarrow H^n(F_t^{-1}(z),\C)_1 \to H^n(T(F_t^{-1}(z)),\C)
\to H^{n-1}_{prim}(F_t)
\end{eqnarray}
The first is the extension to flat sections, the last is a residue map
and is the dual of a tube map.

Consider a form $\omega \in (\Omega_{alg}^{n+1})^{(q)}$ for
$q\in \Z_{\geq (n+1)/d}$ (notation from chapter \ref{s6}).
The restriction of the form $\frac{1}{d}\cdot i_\EE \frac{\omega}{f^q}$
to $F^{-1}_t(z)$ is equal to $\frac{1}{z^{q-1}}\frac{\omega}{\ddd F_t}$, 
because $\ddd F_t\land i_\EE \omega = d\cdot F_t\cdot \omega$.
The section $z\mapsto \left[\frac{1}{z^{q-1}}\frac{\omega}{\ddd F_t}\right]$
of the bundle $\bigcup_{z\in \C^*}H^n(F_t^{-1}(z),\C)_1$ is flat.

Now Varchenko's \cite{Va} description of Steenbrink's Hodge filtration
$F^\bullet$ on $H^\infty_1$ reduces here to the following:
The space $F^{n+1-q}\subset H^\infty_1$ is generated by such flat sections.

In $\P^{n+1}-(F_t^{-1}(z)\cup X_t)$ the set 
$T(F_t^{-1}(z))$ can be deformed to a tubular neighborhood of 
$X_t$ in $\P^n-X_t$. The bundle
$\bigcup_{z\in \C^*}H^n(F_t^{-1}(z),\C)_1$
and the space $H^n(\P^n-X_t)$ glue to a flat bundle on 
$\P^1-\{0\}$. The value in $H^n(\P^n-X_t)$
of a flat section as above is just 
$\left[ \frac{1}{d} \cdot \oooo{i_\EE\frac{\omega}{f^q}}\right]$;
compare the proof of theorem \ref{t6.1}.
This reduces Varchenko's to Griffiths' description and shows (a).

(b) This follows from part (a), theorem \ref{t5.5} and from the following
nontrivial fact: the Frobenius type structure on $T M|_{M_0}$ for a Frobenius
manifold $M$ in theorem \ref{t7.2} corresponds by lemma \ref{t5.1}
to a variation of filtrations and an opposite filtration which are
up to the shift in \eqref{7.9} the variation of Hodge structures and an
opposite filtration on $\bigcup_{t\in M_0}H^\infty_1(F_t)$.

This fact is a consequence of the construction of Frobenius manifolds
in \cite[Theorem 11.1]{He1}. 
\hfill $\qed$

\section{Barannikov--Kontsevich construction}\label{s8}
\setcounter{equation}{0}

\noindent
The Barannikov--Kontsevich construction was initiated in \cite{BK} and 
further developed in \cite{Ba1}\cite{Ba2}. It yields for any
Calabi--Yau manifold a family of formal germs of Frobenius submanifolds.
A central part of it is the construction of a semi-infinite variation
of Hodge structures (defined in \cite{Ba2}). 
This contains the variation of Hodge structures
for complex structure deformations of the given Calabi--Yau manifold.
Therefore it is not surprising that in the case of a Calabi--Yau
hypersurface in $\P^n$ certain submanifolds in the 
Barannikov--Kontsevich construction coincide with germs of the
Frobenius manifolds which one obtains from theorem \ref{t6.1}
with theorem \ref{t5.5}.  
This will be made precise in theorem \ref{t8.1}.
All the results in this chapter are reformulations of results of Barannikov
\cite{Ba1}\cite{Ba2}.

Let $X$ be a compact Calabi--Yau manifold of dimension $n-1$.
Let us fix a holomorphic volume form $\Omega$.
Barannikov \cite{Ba1}\cite{Ba2} constructed a family of formal germs
$(M,0)$ of Frobenius supermanifolds of dimension 
$m:=\dim H^*(X,\C)$. For each of them the tangent space $T_0M$ at 0
is isomorphic to
\begin{eqnarray}\label{8.1}
H^*(X,\bigwedge^*\TT_X) = \bigoplus_{p,q} H^q(X,\bigwedge^p\TT_X)
\end{eqnarray}
with the canonical multiplication and with the $\Z_2$-grading
$(p+q)\mod 2$.
If $t_1,...,t_m$ are flat coordinates centered at 0 with
$\frac{\paa}{\paa t_i}|_0\in H^{q_i}(X,\bigwedge^{p_i}\TT_X)$
then
\begin{eqnarray}\label{8.2}
E = \sum_{i=1}^m \frac{2-p_i-q_i}{2}t_i \frac{\paa}{\paa t_i}
\end{eqnarray}
is an Euler field of the Frobenius supermanifold \cite[5.7]{Ba1}.
Euler field and flat metric $g$ satisfy $\Lie_E(g)=(3-n)g$ and
$g(\frac{\paa}{\paa t_i},\frac{\paa}{\paa t_j})=0$
for $(p_i+p_j,q_i+q_j)\neq (n-1,n-1)$.

This family of Frobenius supermanifolds is parametrized by the set of
{\it opposite filtrations} $W_{\leq \bullet}$ to a Hodge filtration
$F^{\geq \bullet}$ on $H^*(X,\C)$. These filtrations are defined as
follows \cite[ch. 4 and 6]{Ba2}. Both respect the splitting
\begin{eqnarray}\label{8.3}
H^*(X,\C)= H^{even}(X,\C)\oplus H^{odd}(X,\C)
\end{eqnarray}
and are indexed by half integers.
The Hodge filtration is given by
\begin{eqnarray}\label{8.4}
F^{\geq r}:= \bigoplus_{p,q:p-q\geq 2r}H^{p,q}(X) \mbox{ \ for }r\in 
\frac{1}{2}\Z.
\end{eqnarray}
A filtration $W_{\leq \bullet}$ is {\it opposite} if
\begin{eqnarray}\label{8.5}
&& H^*(X,\C) = \bigoplus_{r\in \frac{1}{2}\Z} F^{\geq r}\cap W_{\leq r+1}
\mbox{ \ and }\\
&& (F^{\geq r}\cap W_{\leq r+1},F^{\geq \www r}\cap W_{\leq \www r+1})=0
\mbox{ \ for }r+\www r\neq 0.\label{8.6}
\end{eqnarray}
Here $(,)$ is the Poincar\'e pairing on $H^*(X,\C)$. Let us denote these
Frobenius supermanifolds for a moment by
$M_{Bar}(X,W_{\leq \bullet})$.

From now on we suppose that $n\geq 4$ and that $X$ is a Calabi--Yau 
hypersurface in $\P^n$, that is, $X=\oooo{f^{-1}(0)}\subset \P^n$
where $f\in \C[x_0,...,x_n]$ is homogeneous of degree $n+1$ with 
an isolated singularity at 0. Then the cohomology $H^*(X,\C)$
splits into two orthogonal pieces,
\begin{eqnarray}\label{8.7}
H^*(X,\C) = H^*_{Lef}(X)\oplus H^{n-1}_{prim}(X).
\end{eqnarray}
As in chapter \ref{s5} the second piece 
$H^{n-1}_{prim}(X)\subset H^{n-1}(X)$ is the 
primitive part of the middle cohomology; with $L$ as standard Lefschetz
operator the first piece is
\begin{eqnarray}\label{8.8}
H^*_{Lef}(X) := \bigoplus_{k=0}^{n-1}L^k H^{0}(X,\C).
\end{eqnarray}

If $n-1$ is odd then the splittings \eqref{8.3} and \eqref{8.7} coincide.
Then an opposite filtration $W_{\leq \bullet}$ for $F^{\geq \bullet}$
induces an opposite filtration $\www U_\bullet$ for 
$\www F^\bullet := F^{\bullet-1}H^{n-1}_{prim}(X)$ by
\begin{eqnarray}\label{8.9}
\www U_p := W_{\leq p- \frac{n-1}{2}} \cap H^{n-1}_{prim}(X),
\end{eqnarray}
and this gives a 1--1 correspondence between 
the opposite filtrations $\www U_\bullet$ and the opposite filtrations
$W_{\leq \bullet}$.
If $n-1$ is even then $H^*(X,\C)=H^{even}(X,\C)$. Then formula \eqref{8.9}
gives a 1--1 correspondence between 
the opposite filtrations $\www U_\bullet$ and those opposite filtrations
$W_{\leq \bullet}$ for $F^{\geq\bullet}$ which respect the splitting
\eqref{8.7}.

Let us fix for a moment an opposite filtration $\www U_\bullet$ for 
$\www F^\bullet$
on $H^{n-1}_{prim}(X)$. Theorem \ref{t6.1} and theorem \ref{t5.5}
yield a germ $(M',0)$ of a holomorphic Frobenius manifold of dimension
$m':= \dim H^{n-1}_{prim}(X)$, which is unique up to multiplication of
the metric by a scalar.
Let $E'$ be its Euler field, $\FF'\in \OO_{M',0}$ its potential,
and $e'=\frac{\paa}{\paa t_1}$ its unit field for suitable 
coordinates $t_1,...,t_{m'}$.

The Frobenius manifold can be extended in the following trivial way
to a Frobenius (super)manifold of dimension $m$. Consider 
$\C^{m-m'}$ with coordinates $(\tau_1,...,\tau_{m-m'})$ of degree
$(n-1)\mod 2$. Then the germ
$(M',0)\times (\C^{m-m'},0)$ with potential 
\begin{eqnarray}\label{8.10}
\FF := \FF'+\frac{1}{2}\sum_{i=1}^{m-m'}t_1\tau_i\tau_{m-m'+1-i},
\end{eqnarray}
and Euler field 
\begin{eqnarray}\label{8.11}
E := E' + \sum_{i=1}^{m-m'} \frac{3-n}{2}\tau_i\frac{\paa}{\paa \tau_i}
\end{eqnarray}
is a Frobenius supermanifold for $n-1$ odd and a Frobenius
manifold for $n-1$ even.
We call it $M_{VHS}(X,\www U_\bullet)$. 
It contains the germ $(M',0)\times \{0\}$ as a Frobenius submanifold. 
For $n-1$ odd this is the Frobenius submanifold in \cite[Theorem 8.7.1]{Man2}.
The following theorem is essentially contained in 
\cite[Theorem 6.5]{Ba1}.

\begin{theorem}\label{t8.1}
Let $X\subset \P^n$ be a Calabi--Yau hypersurface with $n\geq 4$.

Consider an opposite filtration $\www U_\bullet$ for 
$\www F^\bullet := F^{\bullet-1}H^{n-1}_{prim}(X)$
on $H^{n-1}_{prim}(X)$ and an 
opposite filtration $W_{\leq\bullet}$ for $F^{\geq\bullet}$ on 
$H^*(X,\C)$ which respects the splitting \eqref{8.7}
if $n-1$ is even. Then
\begin{eqnarray}\label{8.12}
M_{Bar}(X,W_{\leq \bullet}) \cong M_{VHS}(X,\www U_\bullet)
\end{eqnarray}
if and only if \eqref{8.9} holds.
\end{theorem}

{\it Proof.}
For $n-1$ odd $M_{Bar}(X,W_{\leq\bullet})$ is a formal germ of a 
Frobenius supermanifold with odd part corresponding to 
$H^*_{Lef}(X)$ and even part corresponding to 
$H^{n-1}_{prim}(X)$. For $n-1$ even everything is even.
This follows from the discussion at the beginning of this chapter and from
the isomorphisms $\bigwedge^p\TT_X \cong \Omega^{n-1-p}_X$ and 
$H^q(X,\bigwedge^p\TT_X)\cong H^{n-1-p,q}(X)$, which one obtains by
contraction of the holomorphic volume form $\Omega$
on $X$ with holomorphic vector fields.

The condition $n\geq 4$ asserts that $H^{n-2,1}(X)\subset H^{n-1}_{prim}(X)$.
Therefore a miniversal family $X_t$, $t\in M_0$, of complex structure 
deformations of $X_0=X$ is given by a family of homogeneous polynomials
$f_t$, $t\in M_0$, as in chapter \ref{s6} with condition \eqref{6.2}
and $f_0=f$.

By construction of $M_{Bar}(X,W_{\leq\bullet})$ there is a natural inclusion
$(M_0,0)\subset M_{Bar}(X,W_{\leq \bullet})$ of formal germs
(the formal germ $(M_0,0)$ is called $M^{cx}$ in \cite{Ba2}).
By lemma \ref{t4.4} and lemma \ref{t5.1} one obtains on 
$TM_{Bar}(X,W_{\leq \bullet})|_{(M_0,0)}$ a formal germ of a variation
of filtrations with pairing and opposite filtration.

The following fact is crucial: This structure is isomorphic to the 
restriction to the formal germ $(M_0,0)$ of the variation of filtrations
$F^{\geq\bullet}(X_t)$ on the bundle $\bigcup_{t\in M_0}H^*(X_t,\C)$
and the opposite filtration $W_{\leq \bullet}$.

This fact is at the heart of Barannikov's construction of semi-infinite
variations of Hodge structures \cite[Theorem 4.2]{Ba2} and Frobenius manifolds.

Now the case $n-1$ odd is easy. The variation of filtrations 
$F^{\geq \bullet}(X_t)$ and the opposite filtration $W_{\leq \bullet}$
on $\bigcup_{t\in M_0}H^{n-1}_{prim}(X_t)$ correspond to the formal germ
of the even Frobenius submanifold. By theorem \ref{t5.5} this germ
must coincide with the germ $(M',0)$ of the Frobenius manifold before
theorem \ref{t8.1} for filtrations $\www U_\bullet$ and $W_{\leq \bullet}$
with \eqref{8.9}.
Because of the degrees of the flat coordinates the extension with the odd
part to a Frobenius supermanifold is rigid and is given by
\eqref{8.10} and \eqref{8.11}.

The case $n-1$ even is more difficult because the $H^2$-generation 
condition does not hold for the variation of filtrations on 
$\bigcup_{t\in M_0}H^*(X_t,\C)$. But a slightly weaker condition holds so
that lemma \ref{t8.2} below applies. It shows that the variation
of filtrations $F^{\geq \bullet}(X_t)$ on
$\bigcup_{t\in M_0}H^*(X_t,\C)$ and {\it any} opposite filtration
$W_{\leq \bullet}$ determine a unique holomorphic germ of a Frobenius
manifold, whether $W_{\leq \bullet}$  respects the splitting 
\eqref{8.7} or not. It remains to see that precisely the filtrations 
$W_{\leq \bullet}$  which respect the splitting \eqref{8.7} lead to
the Frobenius manifolds $M_{VHS}(X,\www U_\bullet)$.

One can go into the proof of lemma \ref{t8.2} and start there with 
vector fields $\delta_1,...,\delta_m$ and an index set 
$I\subset\{1,...,m\}$ such that the isomorphism
$T_0M\to H^*(X,\C)$ maps $\delta_i$ into $H^*_{Lef}(X)$ for $i\in I$ and
into $H^{n-1}_{prim}(X)$ for $i\notin I$. Then one has to check that
in the case of a filtration $W_{\leq \bullet}$  which respects the
splitting \eqref{8.7}
\begin{eqnarray}\label{8.13}
a_{ij}^k=0\mbox{ for }d_i=0,\ j\in I
\end{eqnarray}
holds. 
Going through the induction in the proof of lemma \ref{t8.2} one finds
that \eqref{8.13} holds 
for all $i\notin I\cup \{1\}$, $j\in I$ and that $a_{ij}^k$ is constant
for $i,j\in I$. Then the flat submanifold $\{t\in M\ |\ t_i =0
\mbox{ for }i\in I\}$ is a Frobenius submanifold, and the
whole Frobenius manifold is $M_{VHS}(\www U_\bullet)$.
\hfill $\qed$

One can generalize theorem \ref{t5.5} by relaxing the $H^2$-generation
condition slightly. The next lemma formulates the part of this generalization
which is needed in the proof of theorem \ref{t8.1}.
It was already used in the proof of \cite[Theorem 6.5]{Ba1}.

\begin{lemma}\label{t8.2}
Consider a germ of a Frobenius manifold $((M,0),\circ,e,E,g)$ with a weight
$w\in \N_{\geq 3}$ and all properties in definition \ref{t5.4}
except the $H^2$-generation condition, which is replaced by the 
weaker condition (II)' on the graded algebra
\begin{eqnarray}\label{8.14}
T_0M&=&\bigoplus_{p=0}^{n-2} (T_0M)_p \mbox{ \ \ with \ }\\
(T_0M)_p &:=& \ker (\nnn^g E-(1-p)\id:T_0M\to T_0M). \nonumber
\end{eqnarray}
(II)' $(T_0M)_1$ generates multiplicatively a subspace of $T_0M$ which
contains $\bigoplus_{p<(w-2)/2} (T_0M)_p$.

Define $M_0 := \{t\in M\ |\ E|_t=0\}$.

The germ of a Frobenius manifold is uniquely determined by the
induced Frobenius type structure (see lemma \ref{t4.4}) on $TM|_{M_0}$
together with $e|_0\in T_0M$.
\end{lemma}

{\it Proof.}
Choose flat coordinates $t_1,...,t_m$ on $(M,0)$ with vector fields
$\delta_i:=\frac{\paa}{\paa t_i}$ such that
\begin{eqnarray}\label{8.15}
E = \sum_{i=1}^m (-d_i)t_i\delta_i
\end{eqnarray}
with $d_1=-1$, $\delta_1=e$, $d_2=...=d_{m_0+1}=0$, $d_i>0$ for
$i>m_0+1$.
Then $M_0=\{t\in M\ |\ t_i=0 \mbox { for }d_i\neq 0\}$.
Define matrices $A_i=(a_{ij}^k)\in M(m\times m,\OO_{M,0})$ for 
$i=1,...,m$ by
\begin{eqnarray}\label{8.16}
\delta_i\circ \delta_j = \sum_k a_{ij}^k\cdot \delta_k.
\end{eqnarray}
Then 
\begin{eqnarray}\label{8.17}
A_iA_j &=& A_jA_i,\\
a_{i1}^k &=&\delta_{ik},\label{8.18}\\
a_{ij}^k &=& a_{ji}^k,\label{8.19}
\end{eqnarray}
and the potentiality condition is equivalent to
\begin{eqnarray}\label{8.20}
\delta_iA_j=\delta_jA_i.
\end{eqnarray}
Especially $\delta_1A_j = 0$ for all $i$. Denote for $w\in \Z_{\geq 0}$
\begin{eqnarray}\label{8.21}
\OO(M)_w &:= & \{f\in \OO_{M_0,0}[t_i\ |\ d_i>0]\ |\ Ef=-wf\},\\
M(w)&:=& M(m\times m,\OO(M)_w), \label{8.22}\\
M(>w) &:=& \bigoplus_{k>w} M(k).\label{8.23}
\end{eqnarray}
Then $t_i\in \OO(M)_{d_i} $ for $i\geq 2$. 
The condition $\Lie_E(\circ) =\circ$ together with \eqref{8.15} shows
\begin{eqnarray}\label{8.25}
a_{ij}^k \in \OO(M)_{d_k-d_i-d_j-1}.
\end{eqnarray}
It is not hard to see that the Frobenius type structure on 
$TM|_{M_0}$ and $e|_0\in T_0M$ provide after some choice the matrices
$A_i\mod M(>0)$ with $d_i=0$ and the coefficients $g(\delta_i,\delta_j)\in\C$
of the metric for all $i,j$.
As in lemma \ref{t2.9} one has to recover all the matrices
$A_j$ in order to uniquely determine the Frobenius manifold.
Again this will be done inductively.

{\it Induction hypothesis for $w\in\Z$}: one has determined the matrices
$A_i\mod M(>w)$ for $d_i=0$ and $A_j\mod M(>w-1)$ for $d_j>0$.

{\it Induction step from $w$ to $w+1$}:
It consists of two steps.
\begin{list}{}{}
\item[(i)]
Determine the matrices $A_j\mod M(>w)$ for $d_j>0$.
\item[(ii)] 
Determine the matrices $A_i\mod M(>w+1)$ for $d_i=0$.
\end{list}

(i) The weakened generation condition (II)' together with \eqref{8.25}
shows that one obtains from the matrices $A_i\mod M(>w)$ with $d_i=0$ 
and from the matrices $A_j\mod M(>w-1)$ the matrices $A_k\mod M(>w)$
for $k$ with $k<\frac{w-2}{2}$. Because of \eqref{8.19} the only
unknown coefficients of the matrices $A_k\mod M(>w)$ for 
$k\geq \frac{w-2}{2}$ are those coefficients $a_{kl}^r$ with 
$d_l\geq \frac{w-2}{2}$. Because of \eqref{8.25}, in the case
$d_k+d_l=w-2$ they are constant and determined by $g(\delta_k,\delta_l)$,
in the case $d_k+d_l>w-2$ they vanish.

(ii) Similarly to step (iii) in the proof of lemma \ref{t2.9} one uses
\eqref{8.20} in the form 
\begin{eqnarray}\label{8.26}
\delta_j(A_i\mod M(>w+1)) =
\delta_i(A_j\mod M(>w+1-d_j))
\end{eqnarray}
for $d_i=0, d_j>0$.
\hfill $\qed$

\begin{remark}\label{t8.3}
%
The proof of theorem \ref{t8.1} contains the three statements:
for $n-1$ even any Frobenius manifold $M_{Bar}(X,W_{\leq \bullet})$
is a holomorphic germ of a Frobenius manifold; it is uniquely determined
up to a scalar of the metric by the variation of filtrations
$F^{\geq\bullet}$ and the opposite filtration $W_{\leq \bullet}$
on $\bigcup_{t\in M_0} H^*(X_t,\C)$; it contains
a Frobenius submanifold with tangent space at 0 isomorphic
to $H^{n-1}_{prim}\subset H^*(X_t,\C)\cong T_0M$ precisely if
the filtration $W_{\leq\bullet}$ respects the splitting \eqref{8.7}
\end{remark}

%
%

\end{document}